\documentclass[11pt]{article}
\usepackage{amsmath,amssymb,theorem}
\usepackage{graphicx, enumerate}
\usepackage{subfigure}
\usepackage{psfrag}
\usepackage{placeins}
\usepackage{wrapfig}
\usepackage{booktabs}
\usepackage{chngcntr}
\usepackage{epstopdf}
\counterwithin{figure}{section}
\numberwithin{figure}{section}
\counterwithin{table}{section}
\usepackage{hyperref}
\hypersetup{
  colorlinks,
  citecolor=black,
  filecolor=black,
  linkcolor=black,
  urlcolor=black
}

\newcommand{\mc}{\mathcal}

\newtheorem{thm}{Theorem}[section]
\newtheorem{conj}[thm]{Conjecture}
\newtheorem{cor}[thm]{Corollary}
\newtheorem{prop}[thm]{Proposition}
\theorembodyfont{\rmfamily}
\def\pf{\bigskip\noindent {\bf Proof.}~~}

\def\dfn#1{{\sl #1}}

\def\less{\backslash}

\def\pf{\bigskip\noindent {\bf{Proof.}}~~}

\newcounter{counter}

\def\proofsquare{  \bigskip\hfill\vrule height3pt width6pt depth2pt}
  \allowdisplaybreaks[4]

\begin{document}
\title{Gallai-Ramsey numbers of $C_{10}$ and $C_{12}$}
\author{Hui Lei$^1$, Yongtang Shi$^1$, Zi-Xia Song$^2$\thanks{Corresponding author.  Email address: Zixia.Song@ucf.edu}\ \ and Jingmei Zhang$^2$\\
{ $^1$ Center for Combinatorics and LPMC}\\
{ Nankai University, Tianjin 300071,  China}\\
{ $^2$ Department  of Mathematics}\\
{ University of Central Florida, Orlando, FL32816, USA}\\
}

\maketitle

\begin{abstract}
A {\it Gallai coloring} is a coloring of the edges of a complete graph without rainbow triangles, and a {\it Gallai $k$-coloring} is a Gallai coloring that uses $k$ colors.  Given   an integer $k\ge1$ and graphs $H_1, \ldots, H_k$, the      {\it Gallai-Ramsey number} $GR(H_1,  \ldots, H_k)$ is the least integer $n$ such that every Gallai $k$-coloring of the complete graph $K_n$   contains a monochromatic copy of $H_i$ in color $i$ for some $i \in \{1,  \ldots, k\}$.  When $H = H_1 = \cdots = H_k$, we simply write $GR_k(H)$.
We  continue to study  Gallai-Ramsey numbers of  even cycles   and paths.    For all $n\ge3$ and $k\ge1$, let  $G_i=P_{2i+3}$ be  a path on $2i+3$ vertices for all $i\in\{0,1, \ldots, n-2\}$ and $G_{n-1}\in\{C_{2n}, P_{2n+1}\}$. Let   $ i_j\in\{0,1,\ldots, n-1 \}$  for all $j\in\{1,  \ldots, k\}$ with
 $  i_1\ge i_2\ge\cdots\ge i_k $.  Song recently conjectured that $ GR(G_{i_1},  \ldots, G_{i_k}) = 3+\min\{i_1,n^*-2\}+\sum_{j=1}^k i_j$, where $n^* =n$ when $G_{i_1}\ne P_{2n+1}$ and $n^* =n+1$ when $G_{i_1}= P_{2n+1}$. This conjecture has been verified to be true for $n\in\{3,4\}$ and all $k\ge1$.  In this paper,  we prove that the aforementioned
conjecture holds for $n \in\{5, 6\}$ and all $k \ge1$. Our result implies that    for all $k\ge1$,
$GR_k(C_{2n})= GR_k(P_{2n})= (n-1)k+n+1$ for $n\in\{5,6\}$  and $GR_k(P_{2n+1})= (n-1)k+n+2$ for      $1\le n \le6 $. \end{abstract}
{\it{Keywords}}: Gallai coloring; Gallai-Ramsey number; Rainbow triangle\\
{\it {2010 Mathematics Subject Classification}}: 05C55;  05D10; 05C15
\section{Introduction}

\baselineskip 16pt

 In this paper we consider graphs that are finite, simple and undirected.    Given a graph $G$ and a set $A\subseteq V(G)$,  we use   $|G|$    to denote  the  number
of vertices    of $G$, and  $G[A]$ to denote the  subgraph of $G$ obtained from $G$ by deleting all vertices in $V(G)\less A$.  A graph $H$ is an \dfn{induced subgraph} of $G$ if $H=G[A]$ for some $A\subseteq V(G)$.  We use $P_n$,  $C_n$ and $K_n$ to denote the path,    cycle and  complete graph  on $n$ vertices, respectively.
For any positive integer $k$, we write  $[k]$ for the set $\{1, \ldots, k\}$. We use the convention   ``$A:=$'' to mean that $A$ is defined to be the right-hand side of the relation.

 Given an integer $k \ge 1$ and graphs $H_1,  \ldots, H_k$, the classical Ramsey number $R(H_1,   \ldots, H_k)$   is  the least    integer $n$ such that every $k$-coloring of  the edges of  $K_n$  contains  a monochromatic copy of  $H_i$ in color $i$ for some $i \in [k]$.  Ramsey numbers are notoriously difficult to compute in general. In this paper, we  study Ramsey numbers of graphs in Gallai colorings, where a \dfn{Gallai coloring} is a coloring of the edges of a complete graph without rainbow triangles (that is, a triangle with all its edges colored differently). Gallai colorings naturally arise in several areas including: information theory~\cite{KG}; the study of partially ordered sets, as in Gallai's original paper~\cite{Gallai} (his result   was restated in \cite{Gy} in the terminology of graphs); and the study of perfect graphs~\cite{CEL}. There are now a variety of papers  which consider Ramsey-type problems in Gallai colorings (see, e.g., \cite{chgr, c5c6,GS, exponential, Hall, DylanSong, C9C11, C13C15}).   These works mainly focus on finding various monochromatic subgraphs in such colorings. More information on this topic  can be found in~\cite{FGP, FMO}.

A \dfn{Gallai $k$-coloring} is a Gallai coloring that uses $k$ colors.
 Given an integer $k \ge 1$ and graphs $H_1,  \ldots, H_k$, the   \dfn{Gallai-Ramsey number} $GR(H_1,  \ldots, H_k)$ is the least integer $n$ such that every Gallai $k$-coloring of $K_n$   contains a monochromatic copy of $H_i$ in color $i$ for some $i \in [k]$. When $H = H_1 = \dots = H_k$, we simply write $GR_k(H)$ and  $R_k(H)$.    Clearly, $GR_k(H) \leq R_k(H)$ for all $k\ge1$ and $GR(H_1, H_2) = R(H_1, H_2)$.    In 2010,
Gy\'{a}rf\'{a}s,   S\'{a}rk\"{o}zy,  Seb\H{o} and   Selkow~\cite{exponential} proved   the general behavior of $GR_k(H)$.

\begin{thm} [\cite{exponential}]
Let $H$ be a fixed graph  with no isolated vertices
 and let $k\ge1$ be an integer. Then
$GR_k(H)$ is exponential in $k$ if  $H$ is not bipartite,    linear in $k$ if $H$ is bipartite but  not a star, and constant (does not depend on $k$) when $H$ is a star.			
\end{thm}

It turns out that for some graphs $H$ (e.g., when $H=C_3$),  $GR_k(H)$ behaves nicely, while the order of magnitude  of $R_k(H)$ seems hopelessly difficult to determine.  It is worth noting that  finding exact values of $GR_k(H)$ is  far from trivial, even when $|H|$ is small.
We will utilize the following important structural result of Gallai~\cite{Gallai} on Gallai colorings of complete graphs.
\begin{thm}[\cite{Gallai}]\label{Gallai}
	For any Gallai coloring $c$ of a complete graph $G$ with $|G|\ge2$, $V(G)$ can be partitioned into nonempty sets  $V_1,  \dots, V_p$ with $p\ge2$ so that    at most two colors are used on the edges in $E(G)\less (E(G[V_1])\cup \cdots\cup  E(G[V_p]))$ and only one color is used on the edges between any fixed pair $(V_i, V_j)$ under $c$.
\end{thm}

The partition given in Theorem~\ref{Gallai} is  a \dfn{Gallai-partition} of  the complete graph $G$ under  $c$.  Given a Gallai-partition $V_1,   \dots, V_p$ of the complete graph $G$ under $c$, let $v_i\in V_i$ for all $i\in[p]$ and let $\mathcal{R}:=G[\{v_1,  \dots, v_p\}]$. Then $\mathcal{R}$ is  the \dfn{reduced graph} of $G$ corresponding to the given Gallai-partition under $c$. Clearly,  $\mathcal{R}$ is isomorphic to $K_p$.
By Theorem~\ref{Gallai},  all edges in $\mathcal{R}$ are colored by at most two colors under $c$.  One can see that any monochromatic $H$ in $\mathcal{R}$ under $c$ will result in a monochromatic $H$ in $G$ under $c$. It is not   surprising   that  Gallai-Ramsey numbers $GR_k(H)$ are closely related to  the classical Ramsey numbers $R_2(H)$.  Recently,  Fox,  Grinshpun and  Pach posed the following  conjecture on $GR_k(H)$ when $H$ is a complete graph.

\begin{conj}[\cite{FGP}]\label{Fox} For all integers $k\ge1$ and $t\ge3$,
\[
GR_k(K_t)= \begin{cases}
			(R_2(K_t)-1)^{k/2} + 1 & \text{if } k \text{ is even} \\
			(t-1)  (R_2(K_t)-1)^{(k-1)/2} + 1 & \text{if } k \text{ is odd.}
			\end{cases}
\]
\end{conj}

The first case of Conjecture~\ref{Fox} follows from a result of     Chung and Graham~\cite{chgr} from  1983.  A simpler proof of this case can be found in~\cite{exponential}.    The    case when $t=4$  was recently settled in~\cite{K4}.  Conjecture~\ref{Fox} remains open for all  $t\ge 5$. The next open case, when $t=5$, involves  $R_2(K_5)$.  Angeltveit and   McKay   \cite{K5} recently  proved that $R_2(K_5)\le 48$. It is widely believed that $R_2(K_5)=43$ (see  \cite{K5}).  It is worth noting that Schiermeyer~\cite{GRK5} recently observed that if  $R_2(K_5)=43$, then Conjecture~\ref{Fox} fails  for $K_5$ when $k=3$.  More recently,  Gallai-Ramsey numbers of  odd cycles on at most $15$ vertices    have been completely settled by Fujita and Magnant~\cite{c5c6} for $C_5$, Bruce and Song~\cite{DylanSong} for $C_7$, Bosse and Song~\cite{C9C11} for $C_9$ and $C_{11}$,  and Bosse, Song and Zhang~\cite{C13C15} for $C_{13}$ and $C_{15}$.

\begin{thm}[\cite{DylanSong, C9C11, C13C15}] For all $k \ge 1$ and $n\in\{3,4,5,6,7\}$,
$GR_k(C_{2n+1}) = n\cdot 2^{k} + 1$.
\end{thm}

In this paper, we continue to study   Gallai-Ramsey numbers of   even cycles and paths.
 For   all  $n \ge 3$ and $k\ge1$, let $G_{n-1}\in \{C_{2n}, P_{2n+1}\}$,   $G_i :=P_{2i+3}$ for all  $i \in \{0, 1, \ldots, n-2\}$, and    $ i_j\in\{0,1,\ldots, n-1 \}$  for all $j\in[k]$. We want to determine  the exact values of $GR (G_{i_1},   \ldots, G_{i_k})$. By reordering colors if necessary, we   assume that $i_1\ge   \cdots \ge i_k$.  Let   $n^* :=n$ when $G_{i_1}\ne P_{2n+1}$ and $n^* :=n+1$ when $G_{i_1}= P_{2n+1}$. Song and Zhang~\cite{C8} recently proved that

 \begin{prop}  [\cite{C8}]\label{lower}
For   all  $n \ge 3$ and $k\ge1$, $$GR(G_{i_1},   \ldots, G_{i_k}) \ge 3+\min\{i_1, n^*-2\}+\sum_{j=1}^k i_j.$$
 \end{prop}

In the same paper, Song~\cite{C8}    further  made the following conjecture.

\begin{conj}[\cite{C8}]\label{Song}
For   all  $n \ge 3$ and $k\ge1$,  
$$GR(G_{i_1},   \ldots, G_{i_k}) = 3+\min\{i_1, n^*-2\}+\sum_{j=1}^k i_j.$$

\end{conj}

To completely solve Conjecture~\ref{Song}, one only needs to consider the case  $G_{n-1}=C_{2n}$.

\begin{prop}[\cite{C8}]\label{C2n to P2n+1}
For   all  $n \ge 3$ and $k\ge1$, if Conjecture~\ref{Song} holds for     $G_{n-1}=C_{2n}$, then
 it also holds for    $G_{n-1}=P_{2n+1}$.
\end{prop}

  Let    $M_n$ denote   a matching of size $n$. As observed in \cite{C8},
the truth of Conjecture~\ref{Song} implies that  $GR_k(C_{2n})=GR_k(P_{2n})=GR_k(M_{n})=(n-1)k+n+1$ for all $n\ge3$ and $k\ge1$ and $GR_k(P_{2n+1})=(n-1)k+n+2$ for all $n\ge1$ and $k\ge1$.    
It is worth noting  that Dzido, Nowik and  Szuca~\cite{R_3(10)} proved that $R_3(C_{2n})\ge4n$ for all $n\ge3$. The truth of Conjecture~\ref{Song} implies that  $GR_3(C_{2n})=4n-2< R_3(C_{2n})$ for all $n\ge3$.  Conjecture~\ref{Song} has recently been verified   to be true for $n\in\{3,4\}$ and all $k\ge1$.

\begin{thm}[\cite{C8}]\label{C8}
For     $n \in\{3,4\}$ and all $k\ge1$,  let  $G_i=P_{2i+3}$ for all  $i \in \{0, 1, \ldots, n-2\}$, $G_{n-1}=C_{2n}$, and     $ i_j\in\{0,1,\ldots, n-1 \}$  for all $j\in[k]$ with $i_1\ge   \cdots \ge i_k$. Then
$$GR(G_{i_1},   \ldots, G_{i_k}) = 3+\min\{i_1, n-2\}+\sum_{j=1}^k i_j.$$
\end{thm}

In this paper, we continue to  establish more evidence for    Conjecture~\ref{Song}.  We prove that Conjecture~\ref{Song} holds  for $n\in\{5,6\}$ and all $k\ge1$.

\begin{thm}\label{main}
For     $n \in\{5,6\}$ and all $k\ge1$,  let  $G_i=P_{2i+3}$ for all  $i \in \{0, 1, \ldots, n-2\}$, $G_{n-1}=C_{2n}$, and     $ i_j\in\{0,1,\ldots, n-1 \}$  for all $j\in[k]$ with $i_1\ge   \cdots \ge i_k$. Then
$$GR(G_{i_1},   \ldots, G_{i_k}) = 3+\min\{i_1, n-2\}+\sum_{j=1}^k i_j.$$
\end{thm}

We prove Theorem~\ref{main} in Section~\ref{section2}. 
We believe the method we developed here can be used to  determine the exact values of $GR_k( C_{2n})$ for  all  $n\ge7$.  Applying Theorem~\ref{main} and Proposition~\ref{C2n to P2n+1}, we obtain the following.

\begin{cor}\label{P11}
Let $G_i=P_{2i+3}$ for all $i\in\{0,1,2,3,4,5\}$. For every  integer  $k \ge 1$, let  $ i_j\in\{0,1,2, 3,4,5 \}$ for all $ j \in[k]$. Then
$$GR (G_{i_1},   \ldots, G_{i_k}) =  3+ \max \{i_j: j \in [k]\} +\sum_{j=1}^k i_j.$$
\end{cor}

\begin{cor}
For all  $k\ge1$,
\begin{enumerate}[(a)]
\item $GR_k (P_{2n+1}) =(n-1)k+n+2$ for all  $n \in [6]$.
\item  $ GR_k (C_{2n}) = GR_k ( P_{2n}) =(n-1)k+n+1$ for $n \in \{5, 6\}$.
\end{enumerate}
\end{cor}

Finally, we shall make use of the following results on 2-colored Ramsey numbers of cycles and paths in the  proof of Theorem~\ref{main}.

\begin{thm}[\cite{Rosta}]\label{cycles}
For all   $n\ge3$, $   R_2(C_{2n}) = 3n-1$.
\end{thm}

\begin{thm}[\cite{FLPS}]\label{path-cycle}
 For all integers $n, m$ satisfying   $2n \geq m \geq 3$,  $R(P_m, C_{2n}) = 2n+\lfloor\frac m 2\rfloor-1$.
\end{thm}

\section{Proof of Theorem~\ref{main}}\label{section2}

We are   ready to prove Theorem~\ref{main}.  Let $n\in\{5,6\}$. By  Proposition~\ref{lower}, it suffices to  show that $ GR(G_{i_1}, \ldots, G_{i_k}) \le 3+\min\{i_1, n-2\}+\sum_{j=1}^k i_j$. \medskip

By Theorem \ref{C8} and Proposition \ref{C2n to P2n+1}, we may assume that $i_1=n-1$. Then $3+\min\{i_1, n-2\}+\sum_{j=1}^k i_j=n+1+\sum_{j=1}^k i_j$. Since $|G_{i_1}|=3+\min\{i_1, n-2\}+  i_1=1+n+i_1$, and $ GR(G_{i_1}, G_{i_2})=R(G_{i_1}, G_{i_2})=1+n+i_1+i_2$ by Theorem~\ref{cycles} and   Theorem~\ref{path-cycle}, we may assume $k\ge3$.
Let $N: =\min\{\max\{i_j:j \in [k]\}, n-2\}+\sum_{j=1}^k i_j$. Then $N \ge 2n-3$.
Let  $G$ be a complete graph on $3+N$ vertices and let $c: E(G)\rightarrow [k]$ be any Gallai  coloring of $ G$ such that all the edges of $G$ are colored by at least three colors under $c$.  We next show that $G$ contains a monochromatic copy of $G_{i_j}$ in color $j$ for some $j\in[k]$. Suppose   $G$ contains no monochromatic copy of $G_{i_j}$ in color $j$ for any $ j\in[k]$ under $c$. Such a Gallai $k$-coloring $c$   is called a  \dfn{bad coloring}.  Among all complete graphs on $3+N$ vertices with a   bad coloring,  we choose $G$ with $N$ minimum. \medskip

Consider a Gallai-partition of $G$ with parts $A_1,   \dots, A_{p}$, where  $p\ge2$.  We may assume that
$|A_1| \ge \cdots  \ge |A_p| \ge 1$.  Let $\mc{R}$ be the reduced graph of $G$ with vertices $a_1,   \ldots, a_p$, where $a_i \in A_i$ for all $i\in[p]$. By Theorem \ref{Gallai}, we may assume that  the edges of $\mc{R}$
are colored either red or blue.   Since all the edges of $G$ are colored by at least three colors under $c$, we see that $\mc{R}\ne G$ and so $|A_1|\ge2$.  By   abusing the notation, we use $i_b$ to denote $i_j$ when  the  color $j$ is blue. Similarly, we use $i_r$ (resp. $i_g$) to denote $i_j$ when  the  color $j$ is red (resp. green).    Let
\[
\begin{split}
A_b &:=  \{a_i \in \{a_2, \ldots, a_p\} \mid a_ia_1 \text{ is colored blue in } \mc{R} \}\\
A_r &:= \{a_j \in \{a_2, \ldots, a_p\} \mid a_ja_1 \text{ is colored red in } \mc{R} \}
\end{split}
\]
Then $|A_b|+|A_r|=p-1$.  Let  $B:= \bigcup_{a_i \in A_b} A_i$ and $R:=\bigcup_{a_j \in A_r} A_j$.  Then $\max\{|B|, |R|\} \ne0$ because $p\ge2$. Thus $G$ contains a blue   $P_3$ between $B$ and $A_1$   or a  red $P_3$ between $R$ and $A_1$, and so  $\max\{i_b, i_r\}\ge 1$.  We next prove several claims. \\

\noindent {\bf Claim\refstepcounter{counter}\label{Observation}  \arabic{counter}.}  Let $r\in[k]$ and let  $s_1, \ldots,s_r$  be nonnegative integers with $s_1+ \cdots+s_r\ge1$.  If $i_{j_1}\geq s_1,  \dots, i_{j_r}\geq s_r$   for colors  $ j_1, \dots, j_r\in[k]$,   then for any $S\subseteq V(G)$ with $|S|\ge|G|-(s_1+ \cdots+s_r)$, $G[S]$ must contain a monochromatic copy of $G_{i^*_{j_q}}$ in color $j_q$ for some $j_q\in\{j_1, \ldots,j_r\}$, where   $i^*_{j_q}=i_{j_q}-s_q$.

\pf  Let $i^*_{j_1} :=i_{j_1}-s_1,  \dots,i^*_{j_r} :=i_{j_r}-s_r$,  and $i^*_j :=i_j$ for all $j\in [k]\less\{j_1, \ldots,j_r\}$. Then  $\max \{i^*_\ell: {\ell \in [k]} \} \le i_1$.  Let $N^* :=\min\{\max \{i^*_\ell: {\ell \in [k]} \}, n-2\}+\sum_{\ell=1}^k i^*_\ell$.   Then $N^*\ge 0$ and $ N^*\le N-(s_1+ \cdots+s_r)<N$ because $s_1+\cdots+s_r\ge1$. Since $|S|\ge 3+N-(s_1+ \cdots+s_r)\ge 3+N^*$ and  $G[S]$ does not have a monochromatic copy of $G_{i_j}$ in color $j$ for all $j\in [k]\less\{j_1, \ldots,j_r\}$ under $c$,  by   minimality of $N$, $G[S]$ must contain a monochromatic copy of $G_{i^*_{j_q}}$ in color $j_q$ for some $j_q\in\{j_1,\ldots,j_r\}$.\proofsquare

 \noindent {\bf Claim\refstepcounter{counter}\label{e:Ap}  \arabic{counter}.}   $|A_1|\le n-1$ and so  $G$ does not contain a  monochromatic copy of  a graph on $|A_1|+1\le n$ vertices in color $m$, where    $m\in[k]$ is  a  color   that is neither red nor blue.

\pf  Suppose $|A_1| \ge n$. We first claim that $i_b\ge |B|$ and $i_r\ge |R|$. Suppose $i_b\le |B|-1$ or  $i_r\le |R|-1$. Then we obtain  a blue  $G_{i_b}$ using the edges between $B$ and $A_1$  or a red $G_{i_r}$ using  the edges  between $R$ and $A_1$, a contradiction.
Thus $i_b\ge |B|$ and $i_r\ge |R|$,  as claimed.   Let $i_b^*:=i_b-|B|$ and  $i_r^*:=i_r-|R|$.
Since $|A_1|= |G|-|B|-|R|$, by Claim~\ref{Observation} applied to $i_b\ge |B|$, $i_r\ge |R|$ and $A_1$, $G[A_1]$ must have a blue $G_{i^*_b}$ or a red   $G_{i^*_r}$. But then either $G[A_1\cup B]$ contains   a blue $G_{i_b}$   or  $G[A_1\cup R]$ contains a  red $G_{i_r}$, because $|A_1| \ge 3+\min\{\max\{i_b, i_r\}, n-2\}+i_b^*+i_r^*$, a contradiction.  This proves that $|A_1|\le n-1$. Next,
let  $m\in[k]$ be  a  color   that is neither red nor blue. Suppose $G$ contains a  monochromatic copy of  a graph, say $J$,  on $n$ vertices in color $m$. Then $V(J)\subseteq A_\ell$ for some $\ell\in[p]$. But then $|A_1|\ge|A_\ell|\ge n$, contrary to $|A_1|\le n-1$. \proofsquare

 For  two disjoint sets $U, W\subseteq V(G)$,  we say $U$  is \dfn{blue-complete} (resp.  \dfn{red-complete})   to $W$    if all the edges between $U$ and $W$    are colored  blue (resp.  red) under $c$.  For convenience, we say $u$  is \dfn{blue-complete} (resp.  \dfn{red-complete})   to $W$  when $U=\{u\}$.\\

\noindent {\bf Claim\refstepcounter{counter}\label{e:R}  \arabic{counter}.}   $\min\{|B|, |R|\} \ge 1$, $p\ge3$,  and $B$ is neither red- nor blue-complete to $R$ under $c$.

\pf  Suppose $B=\emptyset$ or $R=\emptyset$. By symmetry, we may assume that $R=\emptyset$. Then $|B|=|G|-|A_1|=3+N-|A_1|\ge n+1+i_b-|A_1|$.   If $i_b \le |A_1|-1$, then $i_b\le n-2$ by Claim~\ref{e:Ap}. But then we obtain a blue   $G_{i_b}$ using the edges between $B$ and $A_1$. Thus    $i_b \ge |A_1|$. Let $i^*_b=i_b-|A_1|$. By Claim \ref{Observation} applied to $i_b \ge |A_1|$ and $B$, $G[B]$ must have a blue $G_{i^*_b}$. Since  $|B|\ge n+1+i_b^*$, we see that     $G$ contains  a blue   $G_{i_b}$, a contradiction.
Hence $R\ne \emptyset$ and so $p\ge3$ for any Gallai-partition of $G$.   It follows that  $B$ is neither red- nor blue-complete to $R$, otherwise  $\{B\cup A_1, R\}$  or  $\{B, R\cup A_1\}$ yields  a  Gallai-partition of $G$ with only two parts.
 \proofsquare

 \noindent {\bf Claim\refstepcounter{counter}\label{e:ij}  \arabic{counter}.}   Let  $m\in[k]$ be   the color   that is neither red nor blue. Then $i_m\le n-4$. In particular, if $i_m\ge1$, then $G$ contains a monochromatic copy of $P_{2i_m+1}$  in color $m$ under $c$.

\pf By   Claim~\ref{e:Ap},    $|A_1|\le n-1$  and    $G$ contains no   monochromatic copy of $P_{|A_1|+1}$  in color $m$ under $c$.
Suppose $i_m\ge1$.    Let $ i^*_m:=i_m-1$.
 By Claim \ref{Observation} applied to $i_m\ge1$ and  $ V(G) $,  $G$ must have a monochromatic copy of $G_{i^*_m}$  in color $m$ under $c$. Since    $n\in\{5,6\}$, $|A_1|\le n-1$ and $G$ contains no   monochromatic copy of $P_{|A_1|+1}$  in color $m$, we see that $i^*_m \le n-5$. Thus $i_m \le n-4$ and    $G$ contains  a monochromatic copy of $P_{2i_m+1}$  in color $m$ under $c$ if $i_m\ge1$.\proofsquare

By Claim~\ref{e:R} and the fact that $|A_1|\ge2$,   $G$ has a red     $P_3$  and a blue $P_3$. Thus $\min\{i_b, i_r\}\ge1$.  
By Claim~\ref{e:ij}, $\max\{i_b, i_r\}=i_1=n-1$.  Then $|G|= 3+(n-2) +\sum_{i=1}^k i_j \ge 1+n+i_b+i_r \ge 2n+1$.
For the remainder of the proof of Theorem~\ref{main},  we choose  $p\ge3$  to be   as large as possible. \\

   \noindent {\bf Claim\refstepcounter{counter}\label{e:R4}  \arabic{counter}.} $\min\{|B|, |R|\} \le n-1$ if $|A_1| \ge n-3$.

\pf Suppose $|A_1| \ge n-3$ but $\min\{|B|, |R|\} \ge n$.  By symmetry, we may assume that $|B| \ge |R| \ge n$.  Let $B:=\{x_1, x_2, \ldots, x_{|B|}\}$ and $R:=\{y_1, y_2, \ldots, y_{|R|}\}$.  Let $H:=(B,R)$ be the complete bipartite graph obtained from $G[B\cup R]$ by deleting all the edges with both ends in $B$ or  in $R$. Then $H$ has no blue  $P_7$ with both ends in $B$ and no red $P_7$ with both ends in $R$, else we obtain a blue $C_{2n}$ or a red $C_{2n}$ because $|A_1|\ge n-3$. We next show that $H$ has no red $K_{3,3}$. \medskip

Suppose $H$ has a  red $K_{3,3}$. We may assume that   $H[\{x_1, x_2, x_3, y_1, y_2, y_3\}]$ is a red $K_{3,3}$ under $c$.    Since $H$ has   no red $P_7$ with both ends in $R$,  $\{y_4, \ldots, y_{|R|}\}$ must be  blue-complete to $\{x_1, x_2, x_3\}$. Thus   $H[\{x_1, x_2,x_3, y_4, y_5\}]$ has a blue $P_5$ with both ends in $\{x_1, x_2, x_3\}$ and  $H[\{x_1, x_2, x_3, y_1, y_2, y_3\}] $  has  a red $P_5$ with both ends in $\{y_1, y_2, y_3\}$.   If   $|A_1|\ge n-2$ or $\min\{i_b, i_r\}\le n-2$, then we obtain   a blue $G_{i_b}$ or a  red $G_{i_r}$, a contradiction.   It follows that $|A_1|=n-3$ and  $i_b=i_r=n-1$.
Thus $|B \cup R|\ge  1+n+i_b+i_r-|A_1|=2n+2$.   If $|R|\ge6$, then $ \{y_4, y_5, y_6\}$ must be red-complete to $ \{x_4, x_5, x_6\}$, else $H$ has a blue $P_7$ with both ends in $B$. But then we obtain a red $C_{2n}$ in $G$. Thus $|R|=5$, $n=5$, and so $|B|\ge7$.   Let $a_1, a^*_1\in A_1$.  For each $j\in\{4,5,6,7\}$ and every $W\subseteq \{x_1, x_2, x_3\}$ with $|W|=2$, no $x_j$ is red-complete to $W$ under $c$, else,  say,  $x_4$ is red-complete to $\{x_1, x_2\}$, then we obtain a red $C_{10}$  with vertices $a_1, y_1, x_1, x_4, x_2, y_2, x_3, y_3, a^*_1, y_4$ in order, a contradiction.   We may assume that $x_4x_1, x_5x_2$ are colored blue.  But then we obtain a blue $C_{10}$   with vertices $a_1, x_4, x_1, y_4, x_3, y_5, x_2, x_5, a^*_1, x_6$ in order, a contradiction. This proves that $H$ has no red $K_{3,3}$.  \medskip

 Let $X:=\{x_1, x_2, \ldots, x_5\}$ and $Y :=\{y_1, y_2, \ldots, y_5\}$.  Let $H_b$ and  $H_r$ be the spanning subgraphs of $H[X\cup Y]$ induced by all the blue    edges  and red edges of $H[X\cup Y]$ under $c$, respectively.  By the Pigeonhole Principle,  there exist  at least  three   vertices, say  $ x_1, x_2, x_3$, in $ X$ such that either $d_{H_b}(x_i) \ge 3$ for all $i\in[3]$ or $d_{H_r}(x_i) \ge 3$ for all $i\in[3]$. Suppose   $d_{H_r}(x_i) \ge 3$ for all $i\in[3]$. We may assume that $x_1$ is red-complete to $\{y_1, y_2, y_3\}$.   Since $|Y|=5$ and $H$ has no red $P_7$ with both ends in $R$, we see that $N_{H_r}(x_1)=N_{H_r}(x_2)=N_{H_r}(x_3)=\{y_1, y_2, y_3\}$. But then  $H[\{x_1, x_2, x_3, y_1, y_2, y_3\}]$  is  a red $K_{3,3}$, contrary to $H$ has no red $K_{3,3}$.
 Thus  $d_{H_b}(x_i) \ge 3$ for all $i\in[3]$.  Since $|Y|=5$, we see that any  two of $x_1, x_2, x_3$  have a common neighbor in $H_b$.  Furthermore, two of $x_1, x_2, x_3$, say $x_1, x_2$, have at least two common neighbors in $H_b$. It can be easily checked that $H$ has a blue $P_5$ with ends in  $\{x_1, x_2, x_3\}$,  and there exist three vertices, say $y_1, y_2, y_3$, in $Y$ such that $y_ix_i$ is blue for all $i\in[3]$ and $ \{x_4,\ldots,  x_{|B|}\}$ is red-complete to  $\{y_1, y_2, y_3\}$. Then  $H$ has a blue $P_5$ with both ends in $\{x_1, x_2, x_3\}$ and a red $P_5$ with both ends in $\{y_1, y_2, y_3\}$.     If  $|A_1|\ge n-2$ or $\min\{i_b, i_r\}\le n-2$, then we obtain   a blue $G_{i_b}$ or a  red $G_{i_r}$, a contradiction.   It follows that $|A_1|= n-3$ and  $i_b=i_r=n-1$.
 Thus $|B \cup R|\ge 1+n+i_b+i_r-|A_1|=2n+2$. Then $|B|\ge n+1$ and so $H[\{x_4, x_5, x_6, y_1, y_2, y_3\}]$ is  a red $K_{3,3}$, contrary to the fact that $H$ has no red $K_{3,3}$.
  \proofsquare

 \noindent {\bf Claim\refstepcounter{counter}\label{e:A1}  \arabic{counter}.}  $|A_1| \ge 3$. 

\pf Suppose $|A_1| = 2$. Then   $G$ has no monochromatic copy of $P_3$   in color $j$ for any $ j\in\{3, \ldots, k\}$ under $c$. By Claim~\ref{e:ij}, $i_3=\cdots=i_k=0$.
We may assume that  $|A_1| = \cdots =|A_t| =2$ and $|A_{t+1}| =\cdots=|A_p|=1$  for some integer $t$ satisfying $p\ge t\ge1$.  Let $A_i=\{a_i, b_i\}$ for all $i\in [t]$.  By reordering if necessary,  each of $A_1, \ldots, A_t$   can be chosen as the largest part in the Gallai-partition $A_1, A_2,  \ldots, A_p$ of $G$.
For all $i\in [t]$, let
\[
\begin{split}
A^i_b &:=  \{a_j \in V(\mc{R}) \mid a_ja_i \text{ is colored blue in } \mc{R} \}\\
A^i_r &:= \{a_j \in V(\mc{R}) \mid a_ja_i \text{ is colored red in } \mc{R} \}.
\end{split}
\]
 Let $B^i:= \bigcup_{a_j \in A^i_b} A_j$ and $R^i:=\bigcup_{a_j \in A^i_r} A_j$. Then $|B^i|+|R^i|=1+(n-2)+i_b+i_r =2n-2+\min\{i_b, i_r\}$. 
 Let
 \begin{align*}
 E_B &:=  \{a_ib_i \mid i\in [t]   \text{ and }  |R^i|<|B^i| \}  \\
  E_R &:=  \{a_ib_i \mid i\in [t]   \text{ and }  |B^i| <|R^i|\}\\
  E_Q &:= \{a_ib_i \mid i\in [t]   \text{ and }  |B^i|=|R^i|\}.
\end{align*}

Let $c^*$ be   obtained from $c$ by recoloring all the edges in $E_B$   blue, all the edges in $E_R$   red and all the edges in $E_Q$ either red or blue. Then all the edges of  $G$ are colored red or blue under $c^*$. Since $|G|=n+1+i_b+i_r= R(G_{i_b}, G_{i_r})$ by  Theorem~\ref{cycles} and  Theorem~\ref{path-cycle},  we see that $G$ must contain a blue  $G_{i_b}$ or a red $G_{i_r}$ under $c^*$.  By symmetry, we may   assume that $G$ has a blue $H:=G_{i_b}$.  Then $H$ contains no edges   of $E_R$ but must contain  at least one edge of $E_B \cup E_Q$, else we obtain a blue $G_{i_b}$ in $G$ under $c$.  We choose $H$ so that $|E(H) \cap (E_B\cup E_Q)|$ is minimal. We may further assume that $a_1b_1 \in E(H)$. Since $|B^1|+|R^1|=2n-2+\min\{i_b, i_r\}$, by the choice of $c^*$, $|B^1|\ge n-1\geq4$ and $|R^1| \le n-1+\lfloor \frac{\min\{i_b, i_r\}}{2} \rfloor\leq7$. So $i_b\ge2$. By Claim \ref{e:R4}, $|R^1| \le 4$ when $n=5$.  
Let  $W:=V(G) \less V(H) $.\medskip

We next claim that  $i_b=n-1$. Suppose  $i_b\le n-2$. Then $H=P_{2i_b+3}$,  $i_r=n-1$,  $|G|=2n+i_b$ and $|W|=2n-3- i_b \ge n-1$. Let $x_1, x_2, \ldots, x_{2i_b+3}$ be the vertices of $H$ in order.  We may assume that $x_\ell x_{\ell+1}=a_1b_1$ for some $\ell\in [2i_b+2]$.
If a   vertex   $w\in W  $ is blue-complete to $\{a_1, b_1\}$, then  we obtain a blue $H':=G_{i_b}$   under $c^*$ with  vertices $  x_1, \ldots, x_{\ell}, w, x_{\ell+1},  \ldots, x_{2i_b+2}$  in order (when $\ell\ne 2i_b+2$) or $  x_1,   x_{2},    \ldots, x_{2i_b+2}, w$  in order (when $\ell= 2i_b+2$)   such that $|E(H') \cap (E_B\cup E_Q)| < |E(H) \cap (E_B\cup E_Q)|$,  contrary to the choice of $H$.   Thus  no vertex in  $  W$ is blue-complete to $\{a_1, b_1\}$  under $c$ and so $W$ must be red-complete to $\{a_1, b_1\}$ under $c$.  This proves that   $W \subseteq R^1 $.
We next claim that $\ell=1$ or $\ell = 2i_b+2$. Suppose $\ell\in\{2, \ldots, 2i_b+1\}$. Then $\{x_1, x_{2i_b+3}\}$ must be red-complete to $\{a_1, b_1\}$, else, we obtain a blue $H':=G_{i_b}$ with  vertices $x_{\ell}, \ldots, x_1, x_{\ell+1}, \ldots, x_{2i_b+3}$ or $x_1, \ldots, x_{\ell}, x_{2i_b+3}, x_{\ell+1}, \ldots, x_{2i_b+2}$ in order  under $c^*$  such that $|E(H') \cap (E_B\cup E_Q)| < |E(H) \cap (E_B\cup E_Q)|$. Thus  $\{x_1, x_{2i_b+3}\}\subseteq  R^1$ and so $W\cup\{x_1, x_{2i_b+3}\}$ is red-complete to $\{a_1, b_1\}$. If $n=5$, then $4\ge |R^1|\ge |W\cup\{x_1, x_{2i_b+3}\}|\ge 6$, a contradiction. Thus $n=6$ and $7 \ge |R^1|\ge |W\cup\{x_1, x_{2i_b+3}\}|\ge 7$.
It follows that $R^1\cap V(H)=\{x_1, x_{2i_b+3}\}$ and thus either $\{x_{\ell-2}, x_{\ell-1}\}$ or $\{x_{\ell+2}, x_{\ell+3}\}$ is blue-complete to  $\{a_1, b_1\}$. In either case, we  obtain a blue $H':=G_{i_b}$    under $c^*$  such that $|E(H') \cap (E_B\cup E_Q)| < |E(H) \cap (E_B\cup E_Q)|$, a contradiction.
This proves that $\ell=1$ or $\ell = 2i_b+2$.
By symmetry, we may assume that $\ell=1$.  Then $x_1x_3$ is colored blue under $c$ because $A_1=\{a_1, b_1\}$.  Similarly, for all $j\in \{3, \ldots, 2i_b+2\}$,
$\{x_j, x_{j+1}\}$ is not blue-complete to $\{a_1, b_1\}$, else we
obtain a blue $H':=G_{i_b}$ with  vertices $ x_1, x_j, \ldots, x_2,  x_{j+1},  \ldots,  x_{2i_b+3}$  in order  under $c^*$  such that $|E(H') \cap (E_B\cup E_Q)| < |E(H) \cap (E_B\cup E_Q)|$.    It follows that  $x_4\in R^1$ and so $|R^1\cap \{x_4, \ldots, x_{2i_b+3}\}|\ge i_b$.
 Then $|R^1|\ge |W|+|R^1\cap \{x_4, \ldots, x_{2i_b+3}\}|\ge 2n-3$,  so $4 \ge |R^1| \ge 7$ (when $n=5$) or $7 \ge |R^1| \ge 9$ (when $n=6$), a contradiction.  This proves that $i_b=n-1$.
\medskip

\begin{figure}[htbp]
 \centering
 \subfigure[][\label{case1}]{
 \hfill\includegraphics[scale=1]{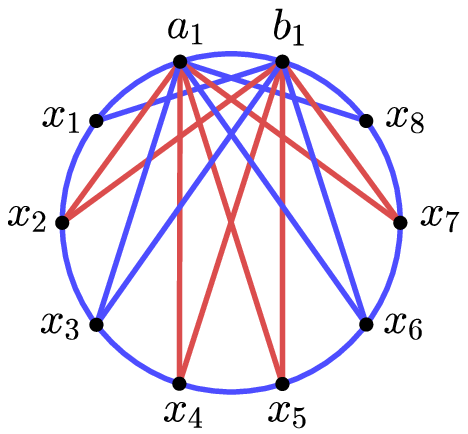}
 }
  \hskip 1cm
 \subfigure[][\label{case2}]{
 \hfill\includegraphics[scale=1]{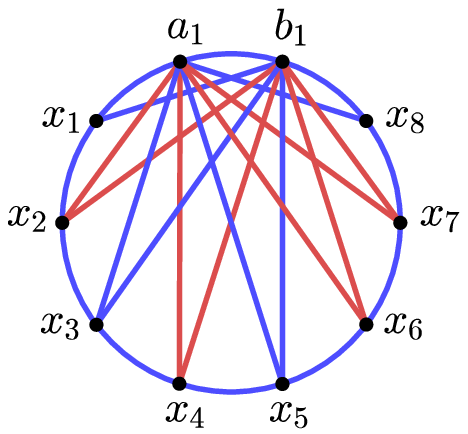}
 }
 \caption{Two cases of $R^1$ when $i_b =4$ and $n=5$.}
 \end{figure}

 Since $i_b=n-1$, we see that  $H=C_{2n}$.   Then $|G|=2n+i_r$ and so $|W|=i_r$.    Let $ a_1, x_1, \ldots,  x_{2n-2}, b_1 $ be the vertices of $H$ in order and let $W:=\{w_1, \ldots,  w_{i_r}\}$. Then  $x_1b_1$ and $a_1x_{2n-2}$ are colored blue under $c$ because $A_1=\{a_1, b_1\}$.
 Suppose  $\{x_j, x_{j+1}\}$ is   blue-complete to $\{a_1, b_1\}$ for some $j\in[2n-3]$. We then obtain a blue $H':=C_{2n}$ with  vertices $a_1, x_1, \ldots, x_j, b_1,$ $ x_{2n-2},\ldots, x_{j+1}$ in order under $c^*$  such that $|E(H') \cap (E_B\cup E_Q)| < |E(H) \cap (E_B\cup E_Q)|$, contrary to the choice of $H$. Thus, for all $j\in[2n-3]$, $\{x_j, x_{j+1}\}$ is not blue-complete to $\{a_1, b_1\}$. Since $\{x_1, x_{2n-2}\}$ is blue-complete to $\{a_1, b_1\}$ under $c$, we see    that $x_2, x_{2n-3}\in R^1$,  and so $4 \ge |R^1\cap V(H)|\ge 4$ (when $n=5$) and $5+\lfloor \frac{i_r}{2} \rfloor \ge |R^1\cap V(H)|\ge 5$ (when $n=6$). 
 Thus,  when $n=5$, we have $R^1=\{x_2, x_4, x_5, x_7\}$ or  $R^1=\{x_2, x_4, x_6, x_7\}$, as depicted in Figure~\ref{case1} and Figure~\ref{case2};
when $n=6$, we have $R^1 \cap V(H)=\{x_2, x_9\} \cup \{x_j: j \in J\}$, where $J\in \{ \{4,6,8\}$, $\{4,6,7\}$, $\{3,4,6,7\}$, $\{3,5,6,7\}$, $\{4,5,6,7\}$, $\{4,6,7,8\}$, $ \{3,5,7,8\}$, $\{3,5,6,8\}$, $\{3,4,5,6,7\}$, $ \{3,4,5,6,8\}$, $\{3,4,5,7,8\} \}$. 

Since $|R^1|\ge n-1$ and $R^1$ is red-complete to $\{a_1, b_1\}$ under $c$, we see that $i_r\ge2$. Let $W':=W \less R^1 \subset B^1$. 
It follows   that  $|W'|=i_r-|R^1 \less V(H)| \ge \lceil \frac{i_r}{2} \rceil \ge 1$. We may assume $W'=\{w_1, \ldots, w_{|W'|}\}$. We claim that $E(H) \cap (E_B \cup E_Q) = \{a_1b_1\}$. Suppose, say     $a_2b_2 \in E(H) \cap (E_B \cup E_Q) $. Since $\{x_1, x_2\}\ne A_i$ and $\{x_{2n-3}, x_{2n-2}\}\ne A_i$ for all $i\in [t]$, we may assume that $a_2=x_j$ and $b_2=x_{j+1}$ for some $j\in \{2, \ldots, 2n-4\}$. Then $x_{j-1}x_{j+1}$ and $x_{j}x_{j+2}$ are colored blue under $c$.  But then we obtain a blue $H':=C_{2n}$  under $c^*$ with vertices $a_1, x_1, \ldots,  x_{j-1}, x_{j+1}, \ldots, x_{2n-2}, b_1, w_1$ in order    such that $|E(H') \cap (E_B \cup E_Q)| < |E(H) \cap (E_B \cup E_Q)|$, contrary to the choice of $H$. Thus $E(H) \cap (E_B \cup E_Q) = \{a_1b_1\}$, as claimed. \medskip

\noindent ($\ast$) Let $w \in W'$. For $j \in \{1,2n-2\}$, if $\{x_j, w\} \neq A_i$  for all $i \in [t]$, then $x_jw$ is colored red. For $j \in \{2,\dots,2n-3\}$, if $\{x_j, w\} \neq A_i$  for all $i \in [t]$ and $x_{j-2}$ or $x_{j+2} \in B^1$, then $x_jw$ is colored red.    

\pf Suppose there are some $j \in [2n-2]$ such that $\{x_j, w\} \neq A_i$  for all $i \in [t]$, and $x_{j-2}$ or $x_{j+2} \in B^1$ if $j \in \{2,\dots,2n-3\}$, but $x_jw$ is colored blue. Then we obtain a blue $C_{2n}$ under $c$ with vertices $a_1, w, x_1, \ldots, x_{2n-2}$ (when $j=1$) or  $a_1, x_1, \ldots, x_{2n-2}, w$ (when $j=2n-2$) in order  if $j\in\{1,2n-2\}$, and  with vertices $b_1, x_{2n-2}, x_{2n-3}, \cdots, x_{j+2}, a_1, w, x_j, \cdots, x_1$ in order (when $x_{j+2} \in B^1$) or $a_1, x_1, \cdots, x_{j-2}, b_1, w, x_j, \cdots, x_{2n-2}$ in order (when $x_{j-2} \in B^1$) if $j \in \{2,\dots,2n-3\}$, a contradiction.  \proofsquare

\noindent ($\ast\ast$) For $j \in [2n-4]$, $x_jx_{j+2}$ is colored red if $\{x_j, x_{j+2}\}\ne A_i$ for all $i\in [t]$.

\pf Suppose $x_jx_{j+2}$ is colored blue for some $j \in [2n-4]$. Then we obtain a blue $C_{2n}$ with vertices $a_1, x_1, \ldots, x_j, x_{j+2}, \ldots, x_{2n-2}, b_1, w$ in order, a contradiction, where $w \in W'$.\proofsquare

First if $n=5$, then $W'=W$.
Let $(\alpha,\beta)\in\{(5,7),(7,6)\}$. Suppose    $R^1=\{x_2, x_4, x_\alpha, x_\beta\}$.
   Since  $\{x_{\alpha-1}, w_j\}\ne A_i$ and $\{x_\alpha, w_j\}\ne A_i$ for all $w_j\in W$ and   $i\in [t]$, $x_{\alpha+1}, x_{\alpha-2} \in B^1$, by ($\ast$),  $\{x_{\alpha-1}, x_\alpha\}$  must be   red-complete to $W$ under $c$.
Then  for any $w_j\in W$, $\{x_{\alpha-2}, w_j\}\ne A_i$ and $\{x_{\alpha+1}, w_j\}\ne A_i$ for all  $i\in [t]$ since $x_{\alpha-1}x_{\alpha-2}$ and $x_\alpha x_{\alpha+1}$ are colored blue  under $c$.  Thus $\{x_{\alpha-2},x_{\alpha+1}\}$ is red-complete to $W$ by ($\ast$). So $\{x_{\alpha-2}, x_{\alpha-1}, x_\alpha, x_{\alpha+1}\}$ is red-complete to $W$ under $c$.
But then we obtain a red $P_9$ under $c$ (when $i_r\le3$) with vertices  $ x_2, a_1,x_{\alpha-1}, b_1, x_\alpha, w_1, x_{\alpha-2}, w_2, x_{\alpha+1}$ in order or a red $C_{10}$ under $c$ (when $i_r=4$) with vertices  $a_1, x_2, b_1, x_{\alpha-1}, w_1, x_{\alpha-2}, w_2, x_{\alpha+1}, w_3, x_{\alpha}$ in order, a contradiction. This proves that $n=6$. By ($\ast$), we may assume $x_1$ is red-complete to $W' \less w_1$ and $x_{10}$ is red-complete to $W' \less w_{|W'|}$ because $|A_1|=2$.
\medskip

\noindent{\bf Case 1.} $|R^1 \cap V(H)|=5$.
 Let $(\alpha,\beta)\in\{(9,8),(7,9)\}$. Suppose  $R^1=\{x_2, x_4, x_6,x_\alpha, x_\beta\}$.
   Since  $\{x_{\alpha-1}, w_j\}\ne A_i$ and $\{x_\alpha, w_j\}\ne A_i$ for all $w_j\in W'$ and   $i\in [t]$,  $x_{\alpha+1}, x_{\alpha-2} \in B^1$, $\{x_{\alpha-1}, x_\alpha\}$  must be   red-complete to $W'$ under $c$ by ($\ast$).  Then for any $w_j\in W'$, $\{x_{\alpha-2}, w_j\}\ne A_i$ and $\{x_{\alpha+1}, w_j\}\ne A_i$ for all  $i\in [t]$ since $x_{\alpha-1}x_{\alpha-2}$ and $x_\alpha x_{\alpha+1}$ are colored blue  under $c$.  Thus $\{x_{\alpha-2},x_{\alpha+1}\}$ is red-complete to $W'$ by ($\ast$). So $\{x_{\alpha-2}, x_{\alpha-1}, x_\alpha, x_{\alpha+1}\}$ is red-complete to $W'$ under $c$. We see that $G$ has a red $P_7$ with vertices $x_{\alpha-1}, w_1, x_{\alpha}, a_1, x_2, b_1, x_4$ in order, and so $i_r \ge 3$ and $|W'| \ge 2$. 
Moreover,  $x_{\alpha-1}x_{\alpha+1}$ and $x_{\alpha-2}x_\alpha$ are colored red by ($\ast\ast$).
Then $G$ has a red $P_{11}$  with vertices  $ x_1, w_2,x_{\alpha-1}, x_{\alpha+1},w_1,x_{\alpha-2},x_\alpha,a_1,x_2,b_1, x_4$ in order under $c$. Thus $i_r=5$ and so $|W'| \ge 3$.  Since $|A_1|=2$ and $x_{\alpha-6} \in B^1$, by ($\ast$), we may assume $x_{\alpha-4}$ is red-complete to $W' \less w_2$. But then we obtain a red $C_{12}$ with vertices  $a_1, x_\alpha, x_{\alpha-2}, w_1,  x_{\alpha-4}, w_3, x_1, w_2, x_{\alpha+1}, x_{\alpha-1},  b_1,  x_2$ in order  under $c$, a contradiction. \medskip

\noindent{\bf Case 2.} $|R^1 \cap V(H)|=6$, then $i_r \ge 3$ and $|W'| \ge 3$.   Let $(\alpha,\beta,\gamma)\in\{(5,2,4),(4,7,5)\}$.
 Suppose $R^1 \cap V(H)=\{x_2, x_3, x_\alpha, x_6, x_7, x_9\}$. 
Since  $\{x_\beta, w_j\}\ne A_i$,  $\{x_3, w_j\}\ne A_i$ and $\{x_6, w_j\}\ne A_i$ for all $w_j\in W'$ and   $i\in [t]$, by ($\ast$),  $\{x_\beta, x_3, x_6\}$  must be   red-complete to $W'$ under $c$. 
By ($\ast\ast$), $x_\gamma$ is red-complete to $\{x_{\gamma-2}, x_{\gamma+2}\}$. 
But then we obtain a red $C_{12}$ under $c$ with vertices $a_1, x_2, x_4, x_6, w_1, x_{10}, w_2, x_1, w_3, x_3, b_1, x_5$ (when  $\alpha=5$) or $a_1, x_3, x_5, x_7, w_1, x_{10}, w_2, x_1, w_3, x_6, b_1, x_4$ (when  $\alpha=4$) in order, a contradiction. 
Let $(\alpha,\beta,\gamma,\delta)\in\{(3,8,5,6), (3,5,7,8),(4,6,8,2)\}$.
 Suppose $R^1 \cap V(H)=V(H) \less \{a_1, b_1, x_1, x_{10}, x_\alpha, x_\beta\}$. Since  $\{x_\gamma, w\}\ne A_i$ and $\{x_\delta, w\}\ne A_i$ for all $w\in W'$ and   $i\in [t]$,  $\{x_\gamma, x_\delta\}$  must be red-complete to $W'$ under $c$ by ($\ast$). Moreover, $x_\gamma x_{\gamma-2}$ and $x_\delta x_{\delta+2}$ are colored red by ($\ast\ast$). Since $|A_1|=2$, there exists at least one of  $x_1, x_{10}, x_\alpha, x_\beta$ is red-complete to $\{w_1, w_2, w_3\}$ by ($\ast$). So we may assume $x_\alpha$ is red-complete to $W' \less w_2$ and $x_\beta$ is red-complete to $\{w_1, w_2, w_3\}$. But then we obtain a red $C_{12}$ with vertices $a_1, x_\gamma, x_{\gamma-2}, w_1, x_{10}, w_2,  x_1, w_3, x_{\delta+2},x_\delta, b_1, x_7 $ in order if $(\alpha,\beta,\gamma,\delta)\in\{(3,8,5,6),(4,6,8,2)\}$ and $a_1, x_7, x_5, w_1, x_3, w_3, x_1, w_2, x_{10}, x_8, b_1, x_6 $ in order if $(\alpha,\beta,\gamma,\delta)=(3,5,7,8)$, a contradiction.
Finally if $R^1 \cap V(H)=\{x_2, x_3, x_5, x_6, x_8, x_9\}$. By ($\ast$), $R^1 \cap V(H)$ is red-complete to $W'$. Then $G$ has a red $P_{11}$ with vertices $x_2, a_1, x_3, b_1, x_5, w_1, x_6, w_2, x_8, w_3, x_9$ in order. Thus $i_r=5$ and so $|W'| \ge 4$. But then we obtain a red $C_{12}$ with vertices $a_1, x_2, w_1, x_3, w_2, x_5, w_3, x_6, w_4, x_8, b_1, x_9$ in order, a contradiction.\medskip

\noindent{\bf Case 3.} $|R^1 \cap V(H)|=7$, then $i_r \ge 4$ and  $|W'| =|W|=i_r$.  Let $(\alpha,\beta)\in\{(6,5), (7,4)\}$. Suppose $R^1 \cap V(H)$=$\{x_2, x_3, x_4, x_5, x_\alpha, x_8, x_9\}$. Since $\{x_3, w_j\}\ne A_i$, $\{x_\beta, w_j\}\ne A_i$ and $\{x_8, w_j\}\ne A_i$ for  all $i\in [t]$ and any $w_j\in W'$,   $\{x_3, x_\beta, x_8\}$  must be red-complete to $W'$ under $c$ by ($\ast$). But then we obtain a red $C_{12}$ with vertices $a_1, x_3, w_1,  x_{10}, w_2, x_1, w_3, x_\beta, w_4, x_8, b_1, x_2$ in order, a contradiction.
Finally if $R^1 \cap V(H)=\{x_2, x_3, x_4, x_5, x_6, x_7, x_9\}$. Since  $\{x_3, w_j\}\ne A_i$ and $\{x_6, w_j\}\ne A_i$ for all   $i\in [t]$ and any $w_j\in W'$,  $\{x_3, x_6\}$  must be red-complete to $W'$ under $c$ by ($\ast$). We may assume $x_8$ is red-complete to $W' \less w_2$ by ($\ast$). But then we obtain a red $C_{12}$ with vertices $a_1, x_3, w_1, x_{10}, w_2, x_1, w_3, x_8, w_4, x_6, b_1, x_2$ in order, a contradiction. This proves   that $|A_1|\ge3$.
\proofsquare

\noindent {\bf Claim\refstepcounter{counter}\label{e:ig}  \arabic{counter}.}  For any $A_i$ with $3 \le |A_i|  \le 4$, $G[A_i]$ has a monochromatic copy of $P_3$ in some color $m \in [k]$ other than red and blue.

\pf  
Suppose there exists a part $A_i$ with $3 \le |A_i| \le 4$ but $G[A_i]$ has no monochromatic  copy  of $P_3$ in any color $m \in [k]$ other than red and blue. We may assume $i=1$. 
Since $GR_k(P_3)=3$, we see that $G[A_1]$ must contain a red or blue $P_3$, say blue.  We may assume $a_i, b_i, c_i$ are the vertices of the blue $P_3$ in order.   
Then  $|A_1|=4$, else $\{b_1\}, \{a_1, c_1\}, A_2, \ldots, A_p$ is a Gallai partition of $G$ with $p+1$ parts.
Let $z_1 \in A_1\less \{a_1, b_1, c_1, \}$. Then $z_1$ is not blue-complete to $\{a_1, c_1\}$, else $\{a_1, c_1\}, \{b_1, z_1\}, A_2, \ldots, A_p$ is a Gallai partition of $G$ with $p+1$ parts. 
Moreover, $b_1z_1$ is not colored blue, else  $\{b_1\}, \{a_1, c_1, z_1\}, A_2, \ldots, A_p$ is a Gallai partition of $G$ with $p+1$ parts. If $b_1z_1$ is colored red, then  $a_1z_1$ and $ c_1z_1$ are colored either red or blue because $G$ has no rainbow triangle. Similarly, $z_1$ is not red-complete to $\{a_1, c_1\}$, else $\{z_1\}, \{a_1, b_1, c_1\}, A_2, \ldots, A_p$ is a Gallai partition of $G$ with $p+1$ parts. Thus, by symmetry, we may assume $a_1z_1$ is colored blue and $c_1z_1$ is colored red, and so $a_1c_1$ is colored blue or red because $G$ has no rainbow triangle. But then $\{a_1\}, \{b_1\}, \{c_1\}, \{z_1\}, A_2, \ldots, A_p$ is a Gallai partition of $G$ with $p+3$ parts, a contradiction. Thus $b_1z_1$ is colored neither red nor blue. But then $a_1z_1$ and $c_1z_1$ must be colored blue because $G[A_1]$ has neither rainbow triangle nor monochromatic $P_3$ in any color $m \in [k]$ other than red and blue, a contradiction. \proofsquare

 For the remainder of the proof of Theorem~\ref{main}, we   assume that $|B|\ge|R|$ . By Claim~\ref{e:R4}, $|R|\le n-1$. Let $\{a_i,b_i,c_i\}\subseteq A_i$ if $|A_i|\geq3$ for any $i\in[p]$. Let $B:=\{x_1, \ldots, x_{|B|}\}$ and $R:=\{y_1, \ldots, y_{|R|}\}$. We next show that \\

\noindent {\bf Claim\refstepcounter{counter}\label{irR}  \arabic{counter}.} $i_r \ge |R|$.

\pf Suppose $i_r \le |R|-1 \le n-2$. Then $i_b=n-1$, $i_r \ge 3$, $|A_1| \le 4$, else we obtain a red $G_{i_r}$ because $R$ is not blue-complete to $B$ and $|A_1| \ge 3$. Moreover, there exist two edges, say $x_1y_1, x_2y_2$, that are colored red, else we obtain a blue $C_{2n}$. Then $G[A_1 \cup R \cup \{x_1, x_2\}]$ has a red $P_9$, it follows that $n=6$, $i_r=4$ and $|R|=5$.  By Claim~\ref{e:ig}, $G[A_1]$ has a monochromatic, say green, copy of $P_3$. By Claim~\ref{e:ij},  $i_g=1$. Then $|A_1 \cup B| =|G|-|R| \ge 7+i_b+i_r+i_g-|R|=12$, and so $G[B]$ has no blue $G_{i_b-|A_1|}$, else we obtain a blue $C_{12}$. Let $i_b^*:=i_b-|A_1|\leq2$, $i_r^*:=i_r-|R|+2=1$, $i_j^*:=i_j \le 2$ for all color $j \in [k]$ other than red and blue.   Let $N^*:=\min\{\max \{i^*_\ell: {\ell \in [k]} \}, 4\}+\sum_{\ell=1}^k i_{\ell}^*$. Observe that $|B|\geq3+N^*$.  By minimality of $N$, $G[B]$ has a red $P_5$ with vertices, say $x_1, \ldots, x_5$, in order. Because there is a red $P_7$ with both ends in $R$ by using edges between $A_1$ and $R$, we see that $R$ is blue-complete to $\{x_1, x_2, x_4, x_5\}$, else $G[A_1 \cup R \cup \{x_1, \ldots, x_5\}]$ has a red $P_{11}$. But then we obtain a blue $C_{12}$ with vertices $a_1, x_1, y_1, x_2, y_2, x_4, y_3, x_5, b_1, x_3, c_1, x_6$ in order, a contradiction. \proofsquare
\medskip

\noindent {\bf Claim\refstepcounter{counter}\label{e:A14}  \arabic{counter}.} $i_b > |A_1|$ and so $|A_1| \le n-2$.

\pf  Suppose $i_b \le |A_1|$.  If $i_b \le |A_1|-1$, then $i_b\le n-2$   by Claim~\ref{e:Ap} and  so $i_r=n-1$. Thus $|B|\ge 2+i_b$ because $|B|+|R|=|G|-|A_1|\ge n+1+i_b+(i_r-|A_1|)  \ge 3+2i_b$.  But then $G$ has a blue   $G_{i_b}$ using edges between $A_1$ and $B$, a contradiction. Thus $i_b=|A_1|$. By Claim~\ref{e:R4} and Claim \ref{irR},   $  |R| \le n-1$ and $i_r \ge |R|$. 
  Observe that $|B|\ge 1+n+i_r-|R|\ge 1+n$. Then     $G[B\cup R]$   has  no blue $P_3$ with both ends in $B$, else we  obtain a blue $G_{i_b}$ in $G$.   Let $i_b^*:=i_b  -|A_1|=0$,   $i_r^*:=i_r-|R|$, and  $i^*_j:=i_j \le n-4$ for all color $j\in [k]$ other than blue and red.  Let $N^*:=\min\{\max \{i^*_\ell: {\ell \in [k]} \}, n-2\}+\sum_{\ell=1}^k i^*_\ell $.   Then $0< N^*<N$.
Suppose first that $|R|\ge2$.  Since $B$ is not red-complete to $R$, we may assume that $y_1x$ is colored blue for some $x\in B$. Note that
$i_r^*  \le n-3$  and  $|B\less x| =3+N-|A_1|-|R|-1\ge 3+N^*$.    By   minimality of $N$,   $G[B \less x]$ must have a red $P_{2i^*_r+3}$ with vertices, say $x_1, \ldots, x_q$, in order, where $q=2i^*_r+3$. Since  $G[B\cup R]$   contains  no blue $P_3$ with both ends in $B$ and $xy_1$ is colored blue, we see that   $y_1$ must be red-complete to  $B \less x$  and   $y_2$ is not blue-complete to $\{x_1, x_q\}$. We may assume that $x_qy_2$ is colored red in $G$. Then $n=6$, $i_r=|R|=5$ and $i_b=|A_1|=3$, else we obtain a red    $G_{i_r}$ using vertices in $V(P_{2i^*_r+3})\cup R\cup A_1$. Let $x' \in B \less \{x, x_1, x_2, x_3\}$. Then $\{x, x'\} \nsubseteq A_i$ and $\{x, x_1\} \nsubseteq A_i$ for all $i \in [p]$ because $yx$ is colored blue and $yx', yx_1$ are colored red, and so $xx'$ and $xx_1$ are colored red, else $G[A_1 \cup B \cup \{y_1\}]$ has a blue $P_9$. But then we obtain a red $C_{12}$ with vertices $a_1, y_1, x', x, x_1, x_2, x_3, y_2, b_1, y_3, c_1, y_4$ in order, a contradiction. Thus $|R|=1$.
  By Claim~\ref{Observation} applied to $i_b=|A_1|$,  $i_r\ge |R|$ and $B$, $G[B ]$ must have a red $P_{2i_r+1}$ with vertices, say $x_1, x_2, \ldots, x_{2i_r+1}$, in order.     Since  $G[B\cup R]$   contains  no blue $P_3$ with both ends in $B$, we may assume that $y_1x_1$ is colored red under $c$. Then $i_r=n-1$, else we obtain a red $G_{i_r}$, a contradiction.  Moreover, $y_1x_{2n-1}$ must be colored blue, else $G$ has a red $C_{2n}$ with vertices $y_1, x_1, \ldots, x_{2n-1}$ in order. Thus $y_1$ is red-complete to $\{x_1, \ldots, x_{2n-2}\}$, and so $\{x_j, x_{2n-1}\} \nsubseteq A_i$ for all $i \in [p]$ and $j \in [2n-2]$.  So $x_{2n-1}x_i$ must be colored red for some $i\in[2n-3]$ because $G[B]$ has no  blue $P_3$. But then we obtain a red $C_{2n}$ with vertices $y_1, x_1, \ldots, x_i, x_{2n-1}, x_{2n-2}, \ldots, x_{i+1}$ in order, a contradiction. 
 This proves that $i_b > |A_1| $, and so $|A_1|\le n-2$.  \proofsquare

By Claim~\ref{e:A1} and Claim~\ref{e:A14}, $3 \le |A_1| \le n-2$. Then by Claim~\ref{e:ig}, $G[A_1]$ has a monochromatic, say green, copy of $P_3$.  By Claim~\ref{e:ij},   $i_g=1$.  \bigskip

\noindent {\bf Claim\refstepcounter{counter}\label{e:A2A3}  \arabic{counter}.} If $|A_1| = 3$, then $|A_2|=3$, $|A_3| \le 2$, and $i_j=0$ for all color $j \in [k] \less [3]$.

\pf  Assume $|A_1|=3$. To prove $|A_2|=3$, we show that $G[B \cup R]$ has  a green $P_3$. Suppose $G[B \cup R]$ has  no green $P_3$. By Claim~\ref{e:A14}, $i_b \ge |A_1|+1=4$. 
Let  $i_g^*:=0$ and  $i^*_j:=i_j$ for all    $j\in [k]$ other than green.  Let $N^*:=\min\{\max \{i^*_\ell: {\ell \in [k]} \}, n-2\}+\sum_{\ell=1}^k i^*_\ell$.   Then $  N^*=N-1$ and $|G\less a_1 | =3+N-1 = 3+N^*$.   But then $G\less a_1$ has no monochromatic copy of $G_{i^*_j}$ in color $j$ for all $j\in[k]$,  contrary to the  minimality of $N$. Thus $G[B \cup R]$ has  a green $P_3$ and so $|A_2|=3$. \medskip

Suppose $|A_3|=3$. For all $i\in [3]$, let 
\[
\begin{split}
A^i_b &:=  \{a_j \in V(\mc{R}) \mid a_ja_i \text{ is colored blue in } \mc{R} \}\\
A^i_r &:= \{a_j \in V(\mc{R}) \mid a_ja_i \text{ is colored red in } \mc{R} \}
\end{split}
\]
 Let $B^i:= \bigcup_{a_j \in A^i_b} A_j$ and $R^i:=\bigcup_{a_j \in A^i_r} A_j$. Since each of $A_1, A_2, A_3$ can be chosen as the largest part in the Gallai-partition $A_1, A_2,\ldots, A_p$ of $G$, by Claim \ref{e:R4}, either $|B^i| \le 5$ or $|R^i| \le 5$ for all $i\in [3]$.
Without loss of generality, we may assume that $A_2$ is blue-complete to $A_1 \cup A_3$.
 Let $X := V(G) \less (A_1 \cup A_2 \cup A_3)=\{v_1, \ldots, v_{|X|}\}$. Then $|X| \ge 1+n+i_b+i_r+i_g-9=2n-8+\min\{i_b, i_r\}$. Suppose $|X \cap B^1| \ge 2$. We may assume  $v_1, v_2 \in X \cap B^1$. Then $G$ has a blue $C_{10}$ with vertices $a_1, v_1, b_1, v_2, c_1, a_2, a_3, b_2, b_3, c_2$ in order and  a blue $P_{11}$ with vertices $a_1, v_1, b_1, v_2, c_1, a_2, a_3, b_2, b_3, c_2, c_3$ in order, and so $n=6$ and $i_b=5$. 
 Moreover, $X \less \{v_1, v_2\} \subseteq R^3$, else, say $v_3$ is blue-complete to $A_3$, then we obtain a blue $C_{12}$ under $c$ with vertices $a_1, v_1, b_1, v_2, c_1, a_2, a_3, v_3, b_3, b_2, c_3, c_2$ in order. Thus $|R^3| \ge |X \less \{v_1, v_2\}|\ge 2+i_r$, and so $i_r \ge 3$, else $G$ has a red $G_{i_r}$  using the edges between $A_3$ and $R^3$. Then there exist at least two vertices in $X \less \{v_1, v_2\}$, say $v_3, v_4$, such that $\{v_3, v_4\}$ is blue-complete to $A_1$, else $G[A_1 \cup A_3 \cup (X \less \{v_1, v_2\})]$ contains a red $G_{i_r}$. Thus $|B^1| \ge |A_2 \cup \{v_1, \ldots, v_4\}|=7$ and so $|R^1| \le 5$.
 Moreover, $\{v_1, v_2\} \subset R^3$, else, say $v_1$ is blue-complete to $A_3$, we then obtain a blue $C_{12}$ under $c$ with vertices $a_1, v_3, b_1, v_4, c_1, a_2, a_3, v_1, b_3, b_2, c_3, c_2$ in order. Then $X \subseteq R^3$ and $|R^3| \ge |X| \ge 4+i_r \ge 7$, and so $|B^3| \le 5$ and  $A_1$ is red-complete to $A_3$.  Furthermore,  $G[B^1\less A_2]$ has no blue $P_3$, else, say $v_1, v_2, v_3$ is such a blue $P_3$ in order, we obtain a blue $C_{12}$ with vertices $a_1, v_1, v_2, v_3, b_1, v_4, c_1, a_2, a_3, b_2, b_3, c_2$ in order.  Therefore for any $U\subseteq B^1 \less A_2$ with $|U|\geq4$, $G[U]$ contains a red $P_3$ because $|A_1|=3$ and $GR_k(P_3)=3$. Since $|R^1| \le 5$ and $A_3 \subseteq R^1$, we may assume $v_1, \ldots, v_{|X|-2} \in B^1 \less A_2$.
Then $G[\{v_1, \ldots, v_4\}]$ must contain a red $P_3$ with vertices, say $v_1, v_2, v_3$, in order. We claim that $X \subset B^1$. Suppose $v_{|X|} \in R^1$. Then $v_{|X|}$ is red-complete to $A_1$ and so $G$ has a red $P_{11}$ with vertices $c_1, v_{|X|}, a_1, a_3, b_1, b_3,  v_1, v_2, v_3, c_3, v_4$ in order, it follows that $i_r=5$.  Thus $|X| \ge 9$, and  $G[\{v_4, \ldots, v_7\}]$ has a red $P_3$ with vertices, say $v_4, v_5, v_6$, in order. But then we obtain a red $C_{12}$ with vertices $a_1, v_{|X|}, b_1, a_3, v_1, v_2, v_3, b_3, v_4, v_5, v_6, c_3$ in order, a contradiction. Thus $X \subset B^1$ as claimed.  Since $|X| \ge 7$, $G[\{v_4, \ldots, v_7\}]$ contains a red $P_3$ with vertices, say $v_4, v_5, v_6$, in order. Then $G$ has a red $P_{11}$ with vertices $a_1, a_3, b_1, b_3, v_1, v_2, v_3, c_3, v_4, v_5, v_6$ in order, and so $i_r=5$, $|X| \ge 9$. Suppose $G[\{v_4, \ldots, v_9\}]$ has no red $P_5$. Then $G[\{v_4, \ldots, v_9\}]$ has at most one part with order three,  say $A_4$, and we may assume $G[A_4]$ has a monochromatic $P_3$ in some color $m$ other than red and blue if $|A_4|=3$ by Claim \ref{e:ig}. 
Let $i^*_r:=1$, $i^*_m:=1$, $i^*_j:=0$ for all color $j\in [k] \less \{m\}$ other than red. Let $N^*:=\min\{i^*_r, 4\}+\sum_{\ell=1}^k i^*_{\ell}=3<N$. Then $G[\{v_4, \ldots, v_9\}]$ has no monochromatic copy of $G_{i_j^*}$ in any color $j \in [k]$, which contradicts to the minimality of $N$. Thus $G[\{v_4, \ldots, v_9\}]$ has a red $P_5$ with vertices, say $v_4, \ldots, v_8$, in order. But then we obtain a red $C_{12}$ with vertices $a_3, v_1, v_2, v_3, b_3, v_4, \ldots, v_8, c_3, v_9$ in order, a contradiction.
 %
Therefore, $|X \cap B^1| \le 1$.  By symmetry, $|X \cap B^3| \le 1$. Let $w \in X \cap B^1$ and $w' \in X \cap B^3$.
 Then $A_1 \cup A_3$ is red-complete to $X \less \{w, w'\}$.  It follows that $n=5$ and $|X \cap B^1|=|X \cap B^3|=1$, else $G[A_1 \cup A_3 \cup (X \less \{w, w'\})]$ has a red $G_{i_r}$ because $|X| \ge 2n-8+\min\{i_b, i_r\}$ and $i_b \ge 4$, a contradiction. But then we obtain a blue $C_{10}$ with vertices $a_2, a_1, w, b_1, b_2, a_3, w', b_3, c_2, c_3$ in order, a contradiction. This proves that $|A_3| \le 2$, and then both $G[A_1]$ and $G[A_2]$ have a green $P_3$, so $i_j=0$ for all color $j \in [k]$ other than  red, blue and green by Claim~\ref{e:ij}.
\proofsquare


\noindent {\bf Claim\refstepcounter{counter}\label{n6R}  \arabic{counter}.}  If $i_b=|A_1|+1$, then $|R|\leq 2$.
   
\pf Suppose $i_b=|A_1|+1$ but $|R|\geq3$. By Claim~\ref{irR}, $i_r \ge |R|$, it follows that $|B| \ge 1+n+i_b+i_r+i_g-|A_1|-|R| \ge 3+n$. Thus $G[B \cup R]$ has no blue $P_5$ with both ends in $B$, else we obtain a blue $G_{i_b}$. Let $i_b^*:=i_b-|A_1| = 1$, $i_r^*:=i_r-|R|+1$ (when $n=5$) or $i_r^*:=\max\{i_r-|R|+1,2\}$ (when $n=6$), $i_j^*:=i_j$ for all $j \in [k]$ other than red and blue. Let $N^*:=\min\{\max \{i^*_\ell: {\ell \in [k]} \}, n-2\}+\sum_{\ell=1}^k i_{\ell}^*$. Then $0< N< N^*$. Observe that $|B|\geq3+N^*$. By minimality of $N$, $G[B]$ has  a red $G_{i_r^*}$ with vertices, say $x_1, \ldots, x_q$, in order, where $q=2i_r^*+3$.  If $R$ is blue-complete to $\{x_1, x_q\}$, then $R$ is red-complete to $B \less \{x_1, x_q\}$ because $G[B \cup R]$ has no blue $P_5$ with both ends in $B$. But then $G[A_1 \cup R \cup \{x_2, \ldots, x_{q-1}\}]$ has a red $G_{i_r}$, a contradiction. Thus $R$ is not blue-complete to $\{x_1, x_q\}$, and so we may assume $y_1x_1$ is colored red. Then $i_r=n-1$ and $R\less\{y_1\}$ is blue-complete to $ \{x_{q-2}, x_q\}$, else $G[A_1 \cup R \cup \{x_1, \ldots, x_q\}]$ has a red $G_{i_r}$. So $R\less\{y_1\}$ is  red-complete to $B \less \{x_{q-2}, x_q\}$ because $G[B \cup R]$ has no blue $P_5$ with both ends in $B$. But then $G[A_1 \cup R \cup \{x_2, \ldots, x_{q-1}\}]$ has a red $G_{i_r}$, a contradiction. \proofsquare

\noindent {\bf Claim\refstepcounter{counter}\label{e:ibn}  \arabic{counter}.} $i_b = n-1$.

\pf Suppose $i_b \le n-2$. Then $i_r=5$. By Claim \ref{e:A1} and Claim \ref{e:A14}, $|A_1| \ge 3$ and $i_b >|A_1|$, it follows that $n=6$, $i_b=4$ and $|A_1|=3$. By Claim \ref{e:A2A3},  $|A_2|=3$, $|A_3| \le 2$, $i_j=0$ for all color $j \in [k] \less [3]$. By Claim \ref{n6R}, $|R| \le 2$ and so $A_2 \subset B$.  It follows that $|B|= 7+i_b+i_r+i_g-|A_1 \cup R| =14-|R| \ge 12$.  Then $G[B \cup R]$ has no blue $P_5$ with both ends in $B$, else $G$ has a blue $P_{11}$ because $|A_1|=3$.
Thus there exists a set $W$ such that $ (B \cup R) \less (A_2 \cup W)$ is red-complete to $A_2$, where $W \subset (B \cup R) \less A_2$ with $|W|\le 1$. Let $i_b^*:=i_b-|A_1| = 1$, $i_r^*:=2$, $i_j^*:=0$ for all $j \in [k]$ other than red and blue. Let $N^*:=\min\{\max \{i^*_\ell: {\ell \in [k]} \}, 4\}+\sum_{\ell=1}^k i_{\ell}^*=5$. Then $N^* < N$. Observe that $|B \less (A_2 \cup W)| = 11-|R|-|W| \ge 3+N^*$. By minimality of $N$, $G[B \less (A_2 \cup W)]$ must contain a red $G_{i_r^*}=P_7$. But then $G[(B \cup R) \less W]$ has a red $C_{12}$, a contradiction. Thus $i_b=n-1$.
\proofsquare

\noindent {\bf Claim\refstepcounter{counter}\label{e:A2}  \arabic{counter}.} $|A_1| = n-2$.

\pf By Claim~\ref{e:A14}, $|A_1| \le n-2$.  Suppose $|A_1| \le n-3$. By Claim~\ref{e:A1}, $n=6$ and $|A_1|=3$. By Claim~\ref{e:ibn}, $i_b =5$. By Claim \ref{e:A2A3}, $|A_2|=3$, $|A_3| \le 2$ and $i_j=0$ for all color $j \in [k]\less [3]$. By Claim \ref{irR}, $i_r \ge |R|$. Then $|B| = 7+i_b+i_r+i_g-|A_1|-|R|\ge 10$, and so $G[B \cup R]$ has neither blue $P_7$ nor blue $P_5 \cup P_3$ with all ends in $B$ else we obtain a blue $C_{12}$. \medskip

Suppose $|R| \le 2$. Then $A_2 \subset B$ and there exists a set $W \subset (B  \cup R) \less A_2$ with $|W| \le 3$ such that $W$ is blue-complete to $A_2$ and $(B  \cup R) \less (A_2 \cup W)$ is red-complete to $A_2$. Since $|B \less (A_2 \cup W)| \ge 4$, we see that there is a red $P_7$ using edges between $A_2$ and $B \less (A_2 \cup W)$, so $i_r \ge 3$ and $i_r-|R| \ge 1$. Let $i_b^*:=2$ (when $|B \cap W| \le 1$) or $i_b^*:=0$ (when $|B \cap W| \ge 2$), $i_r^*:= \min\{i_r-|R|-1,2\}$, $i_j^*:=0$ for all color $j \in [k]$ other than red and blue,  and $N^*:=\min\{\max \{i^*_\ell: {\ell \in [k]} \}, 4\}+\sum_{\ell=1}^k i_{\ell}^*=\max\{i_b^*, i_r^*\}+i_b^*+i_r^*$. 
Observe that $|B \less (A_2 \cup W)|=7+i_r-|R \cup W|\ge 3+N^*$. By minimality of $N$, $G[B \less (A_2 \cup W)]$ has a red $G_{i_r^*}$ because $G[B]$ has neither blue $P_7$ nor blue $P_5 \cup P_3$ and $|A_3| \le 2$. But then $G[(B\cup R) \less W]$ has a red $G_{i_r}$ because $|(B\cup R) \less W| \ge 7+i_r \ge |G_{i_r}|$ and $A_2$ is red-complete to $(B  \cup R) \less (A_2 \cup W)$, a contradiction.
Therefore, $3 \le |R| \le 5$ and so $i_r \ge 3$. \medskip

We claim that $i_r=5$. Suppose $3 \le i_r \le 4$. Let $i_b^*:=2$, $i_r^*:=2$, $i_j^*:=i_j$ for all color $j \in [k]$ other than red and blue,  and $N^*:=\min\{\max \{i^*_\ell: {\ell \in [k]} \}, 4\}+\sum_{\ell=1}^k i_{\ell}^*=7$. Observe that $|B| \ge 10=3+N^*$. Since $G[B]$ has no blue $P_7$, by minimality of $N$, $G[B]$ has a red $P_7$ with vertices, say $x_1, \ldots, x_7$, in order. Then $R$ is blue-complete to $\{x_1, \ldots, x_7\} \less x_4$, else $G[A_1 \cup R \cup \{x_1, \ldots, x_7\}]$ has a red $G_{i_r}$. But then $G[B \cup R]$ has a blue $P_7$ with vertices $x_1, y_1, x_2, y_2, x_3, y_3, x_5$ in order, a contradiction. Thus $i_r=5$ and so $|G|=18$, $|B|=15-|R|$. \medskip

If $|R|=3$. First suppose $A_2 \subseteq R$. Since $R$ is not red-complete to $B$, we may assume that $A_2$ is blue-complete to $x_1$. Let $i_b^*:=2$,   $i_r^*:=3$, $i_j^*:=0$ for all color $j \in [k]$ other than red and blue,  and $N^*:=\min\{\max \{i^*_\ell: {\ell \in [k]} \}, 4\}+\sum_{\ell=1}^k i_{\ell}^*=8$. Observe that $|B \less x_1|=11=3+N^*$. By minimality of $N$, $G[B \less x_1]$ has a red $P_9$ with vertices, say $x_2, \ldots, x_{10}$, in order. We claim that $A_2$ is blue-complete to $\{x_2, x_{10}\}$, else, say $x_2$ is red-complete to $A_2$. Then $A_2$ is blue-complete to $\{x_8, x_{10}\}$, else $G[A_1 \cup A_2 \cup \{x_2, \ldots, x_{10}\}]$ has a red $C_{12}$. Thus $A_2$ is red-complete to $B \less \{x_1, x_8, x_{10}\}$ because $G[B \cup R]$ has no blue $P_7$ with both ends in $B$. But then we obtain a red $C_{12}$ with vertices $a_1, a_2, x_3, \ldots, x_9, b_2, b_1, c_2$ in order, a contradiction. Thus,  $A_2$ is blue-complete to $\{x_1, x_2, x_{10}\}$, and so $A_2$ is red-complete to $B \less \{x_1, x_2, x_{10}\}$ because $G[B \cup R]$ has no blue $P_7$ with both ends in $B$. But then we obtain a red $C_{12}$ with vertices $a_1, a_2, x_3, \ldots, x_9, b_2, b_1, c_2$ in order, a contradiction.  This proves that $A_2 \subset B$. Then there exists a set $W \subset (B  \cup R) \less A_2$ with $|W \cap B| \le 3$ such that $W$ is blue-complete to $A_2$ and $(B  \cup R) \less (A_2 \cup W)$ is red-complete to $A_2$. Then  $|W| \le 3$ and $|W \cap B| \le 3$ or $|W| = 4$ and $|W \cap B| = 1$ because $G[B \cup R]$ has no blue $P_7$ with both ends in $B$. Let
\begin{align*}
& i_b^*:=2-|W|,\ i_r^*:=2 \,\ \text{when}\ |W| \in \{0, 1\} \\
& i_b^*:=0, \ i_r^*:=2 \,\ \text{when}\ |W| \ge 2  \ \text{and }\ |W \cap B| \le 2 \\
& i_b^*:=0, \ i_r^*:=1 \,\ \text{when}\ |W|=|W \cap B| =3,
\end{align*} 
$i_j^*:=0$ for all color $j \in [k]$ other than red and blue,  and $N^*:=\min\{\max \{i^*_\ell: {\ell \in [k]} \}, 4\}+\sum_{\ell=1}^k i_{\ell}^*=2i_r^*+i_b^*$. Observe that $|B \less (A_2 \cup W)| \ge 3+N^*$. By minimality of $N$, $G[B \less (A_2 \cup W)]$ has a red $G_{i_r^*}$ because $G[B \cup R]$ has neither blue $P_7$ nor blue $P_5 \cup P_3$ with all ends in $B$ and $|A_3| \le 2$. If $|W| \le 3$ and $|W \cap B| \le 2$, then $G[(B \cup R) \less W]$ has a red $C_{12}$ because $|(B \cup R) \less W|\geq12$ and $A_2$ is red-complete to $(B\cup R) \less (A_2 \cup W)$. Thus $|W|=|W \cap B|=3$ or $|W|=4$ and $|W \cap B|=1$. For the former case, $G[B \less (A_2 \cup W)]$ has a red $P_5$ with vertices, say $x_1, \ldots, x_5$, in order. Let $W:=\{w_1, w_2, w_3\} \subset B$. Then $A_2$ is blue-complete to $W$ and red-complete to $\{x_1, \ldots, x_5\}$, and so $W$ is red-complete to $\{x_1, \ldots, x_5\}$ because $G[B]$ has no blue $P_7$. But then we obtain a red $C_{12}$ with vertices $a_2, x_1, w_1, x_2, w_2, x_3, w_3, x_4, b_2, x_5, c_2, x_6$ in order, where $x_6 \in B \less (A_2 \cup W \cup \{x_1, \ldots, x_5\})$, a contradiction. For the latter case, $G[B \less (A_2 \cup W)]$ has a red $P_7$ with vertices, say $x_1, \ldots, x_7$, in order. Let $W \cap B:=\{w\}$. Then $w$ is red-complete to $\{x_1, \ldots, x_7\}$ because $G[B]$ has no blue $P_7$. But then we obtain a red $C_{12}$ with vertices $a_2, x_1, w, x_2, \ldots, x_6, b_2, x_7, c_2, x_8$ in order, where $x_8 \in B \less (A_2 \cup W \cup \{x_1, \ldots, x_7\})$, a contradiction.
This proves that $|R| \in \{4, 5\}$.  First we claim that $G[E(B, R)]$ has no blue $P_5$ with both ends in $B$. Suppose there is a blue $H:=P_5$ with vertices, say $x_1, y_1, x_2, y_2, x_3$, in order.  Then  $G[(B \cup R) \less V(H)]$ has no blue $P_3$ with both ends in $B$. Let $i_b^*:=0$, $i_r^*:=i_r-|R|+1$, $i_j^*:=i_j$ for all color $j \in [k]$ other than red and blue,  and $N^*:=\min\{\max \{i^*_\ell: {\ell \in [k]} \}, 4\}+\sum_{\ell=1}^k i_{\ell}^*=3+2(i_r-|R|)$. Observe that $|B \less \{x_1, x_2, x_3\}| = 7+i_r-|R| \ge 3+N^*$ since $|R| \in \{4, 5\}$. By minimality of $N$, $G[B \less \{x_1, x_2, x_3\}]$ has a red $G_{i_r^*}$ with vertices, say $x_4, \ldots, x_{q}$, in order, where $q=2i_r^*+6$. Then $y_3$ is not blue-complete to $\{x_4, x_q\}$ because $G[(B \cup R) \less V(H)]$ has no blue $P_3$ with both ends in $B$. We may assume $x_4y_3$ is colored red. Then $R \less \{y_1,y_2,y_3\}$ is blue-complete to $x_8$, else we obtain a red $C_{12}$ with vertices $a_1, y_3, x_4, \ldots, x_8, y_4, b_1, y_1, c_1, y_2$ in order, a contradiction. Since $G[(B \cup R) \less V(H)]$ has no blue $P_3$ with both ends in $B$, we see that $R \less \{y_1,y_2,y_3\}$ is red-complete to $\{x_4, \ldots, x_q\} \less \{x_8\}$. But then we obtain a red $C_{12}$ with vertices $a_1, y_3, x_4, \ldots, x_{10}, y_4, b_1, y_1$ (when $|R|=4$), or $a_1, y_3, x_4, x_5, x_6, y_4, x_7, y_5, b_1, y_1, c_1, y_2$ (when $|R|=5$) in order, a contradiction. Thus, $G[E(B, R)]$ has no blue $P_5$ with both ends in $B$.
%
Let $i_b^*:=2$, $i_r^*:=2$, $i_j^*:=i_j$ for all color $j \in [k]$ other than red and blue,  and $N^*:=\min\{\max \{i^*_\ell: {\ell \in [k]} \}, 4\}+\sum_{\ell=1}^k i_{\ell}^*=7$. Observe that $|B| \ge 10=3+N^*$. By minimality of $N$, $G[B]$ has a red $P_7$ with vertices, say $x_1, \ldots, x_7$, in order. We claim that $x_1$ is blue-complete to $R$. Suppose $x_1y_1$ is colored red. Then $R \less y_1$ is blue-complete to $\{x_5, x_7\}$, else $G[A_1 \cup R \cup \{x_1, \ldots, x_7\}]$ has a red $C_{12}$. Thus $R \less y_1$ is red-complete to $B \less \{x_5, x_7\}$ because $G[E(B, R)]$ has no blue $P_5$ with both ends in $B$. But then we obtain a red $C_{12}$ with vertices $a_1, y_2, x_2, \ldots, x_6, y_3, b_1, y_4, c_1, y_1$ in order, a contradiction. Therefore, $x_1$ is blue-complete to $R$. By symmetry, $x_7$ is blue-complete to $R$. Then $R$ is red-complete to $B \less \{x_1, x_7\}$ because $G[E(B, R)]$ has no blue $P_5$ with both ends in $B$. But then we obtain a red $C_{12}$ with vertices $a_1, y_2, x_2, \ldots, x_6, y_3, b_1, y_4, c_1, y_1$ in order, a contradiction.  This proves that $|A_1| = n-2$. \proofsquare

By Claim~\ref{e:ibn}, Claim \ref{e:A2} and Claim \ref{irR},   $i_b=n-1$,  $|A_1|=n-2$,   $ i_r\ge |R|$. By Claim \ref{n6R}, $|R|\le 2$. Then $|B|\ge 3+n+i_r-|R|\ge 3+n$, and so $G[B\cup R]$ has no blue $P_5$ with both ends in $B$.
 \medskip 

\noindent {\bf Claim\refstepcounter{counter}\label{e:ir=4}  \arabic{counter}.} $i_r=n-1$.

\pf Suppose $i_r\le n-2$. By Claim~\ref{e:R}, $B$ is not blue-complete to $R$. Let $x\in B$ and $y\in R$ such that $xy$ is colored red.
Let $i_b^*:=i_b  -|A_1|=1$ and $i_r^*:=i_r-|R|\le n-3$,  $i^*_j:=i_j\le n-4$ for all color  $j\in [k]$ other than red and blue.  Let $N^*:=\min\{\max \{i^*_\ell: {\ell \in [k]} \}, n-2\}+\sum_{\ell=1}^k i^*_\ell$.   Then $0< N^*<N$ and $|B\less x| =3+N-|A_1|-|R|-1\ge 3+N^*$.    By   minimality of $N$,   $G[B \less x]$ must have a red $P_{2i^*_r+3}$ with vertices, say $x_1, x_2, \ldots, x_{2i^*_r+3}$, in order. 
Then $\{x_1, x_{2i^*_r+3}\}$ must be blue-complete to $\{x,y\}$ and $xx_2$ must be colored blue under $c$,  else we obtain a red $P_{2i_r+3}$ using vertices in $V(P_{2i^*_r+3})\cup\{x,y\}$ or  in $V(P_{2i^*_r+3}\less x_1)\cup\{x,y\}\cup A_1$.  But then  $G[B\cup R]$ has a blue $P_5$ with vertices $x_2, x, x_1, y, x_{2i^*_r+3}$ in order, a contradiction.
\proofsquare

 Let   $A_1:=\{a_1, b_1, c_1\}$ (when $n=5$) or $A_1:=\{a_1, b_1, c_1, z_1\}$ (when $n=6$). 
 By Claim~\ref{e:ig}, $G[A_1]$ has a monochromatic, say green, copy of $P_3$. By Claim~\ref{e:ij},   $i_g=1$.  We next show  that   $|A_2| \ge 3$. Suppose $|A_2| \le 2$.
 Then by Claim~\ref{e:A2A3}, $|A_1|=4$ and so $n=6$.  Let $i_b^*:=i_b-|A_1|$, $i_r^*:= i_r-|R|+1$, $i_g^*:=i_g-1=0$ and  $i^*_j:=i_j$ for all    $j\in [k]$ other than red, blue and green.  Let $N^*:=\min\{\max \{i^*_\ell: {\ell \in [k]} \}, 4\}+\sum_{\ell=1}^k i^*_\ell$.  Then $ 0< N^*<N$ and $|B| =|G|-|A_1|-|R| = 3+N^*$. By minimality of $N$, $G[B]$ must contain a red $G_{i^*_r}$. It follows that $|R|=2$ and $G_{i^*_r}= P_{11}$. Let $x_1, x_2, \ldots, x_{11}$ be the vertices of the red $P_{11}$ in order. If $R$ is blue-complete to $\{x_1, x_{11}\}$, then $R$ is red-complete to $B \less \{x_1, x_{11}\}$ because $G[B\cup R]$ has no blue $P_5$ with both ends in $B$.  But then $G$ has a red $C_{12}$ with vertices $a_1, y_1, x_2, \ldots, x_{10}, y_2$ in order, a contradiction. Thus, $R$ is not blue-complete to $\{x_1, x_{11}\}$ and we may assume $x_1y_1$ is colored red. Then $x_{11}y_1$ and $x_9y_2$ are colored blue, else $G[\{x_1, \ldots, x_{11}\} \cup R \cup A_1]$ has a red $C_{12}$. If $x_{11}y_2$ is colored red, then $x_1y_2$ and $x_3y_1$ are colored blue by the same reasoning. But then we obtain a blue $C_{12}$ with vertices $a_1, x_1, y_2, x_9, b_1, x_3, y_1, x_{11}, c_1, x_2, z_1, x_4$ in order, a contradiction. Thus $x_{11}y_2$ is colored blue. Then $y_1$ is red-complete to $B \less \{x_9, x_{11}\}$, else, say  $y_1w$  is colored blue with $w \in B \less \{x_9, x_{11}\}$,  then $G[B \cup R]$ has  a blue $P_5$ with vertices $w, y_1, x_{11}, y_2, x_9$ in order. It follows that $\{x_{11}, w\} \nsubseteq A_j$ for all $j \in [q]$, where $w \in B \less \{x_9, x_{11}\}$. Moreover, $x_2y_2$ is colored blue, else $G$ has a red $C_{12}$ with vertices $a_1, y_2, x_2, \ldots, x_{10}, y_1$ in order, a contradiction. Thus, $G[B \less \{x_2, x_9\}]$ has no blue $P_3$, else $G[A_1 \cup B \cup \{y_2\}]$ has a blue $C_{12}$. Therefore, $x_ix_{11}$ is colored red for some $i \in \{3, \ldots, 7\}$. But then we obtain a red $C_{12}$ with vertices $y_1, x_1, \ldots, x_i, x_{11}, x_{10}, \ldots, x_{i+1}$ in order, a contradiction.  Thus  $3 \le |A_2| \le n-2$ and $A_2 \subset B$ because $|R| \le 2$.\medskip

 Since $G[B\cup R]$ has no blue $P_5$ with both ends in $B$,  there exists at most one  vertex, say $w \in B\cup R$, such that  $(B\cup R)\less (A_2 \cup\{w\})$ is  red-complete to $A_2$.   Suppose  $3 \le |A_3| \le n-2$.  Then $n=6$ by Claim~\ref{e:A2A3}, $A_3\subseteq B$ and  $A_3$ must be red-complete to $A_2$.    Since $G[B\cup R]$ has no blue $P_5$ with both ends in $B$, there exists at most one vertex, say $  w'\in B\cup R$, such that    $(B\cup R)\less (A_3 \cup\{w'\})$ is  red-complete to $A_3$.  But then $G[(B \cup R)\less \{w,w'\} ]$ has  a red $C_{12}$, a contradiction.
 Thus $|A_3|\le2$ and so $G[B\less A_2]$ has no monochromatic  copy of $P_3$ in color $j$ for all $j\in[k]$ other than red and blue.
  Let $i_b^*:=1$, $i_r^*:=n-1-|A_2|$,  and  $i^*_j:=0$ for all colors  $j\in [k]$ other than red and blue.  Let $N^*:=\min\{\max \{i^*_\ell: {\ell \in [k]} \},n-2\}+\sum_{\ell=1}^k i^*_\ell=2i_r^*+1=2n-1-2|A_2|$.   Then $0< N^*<N$ and $|B \less (A_2\cup \{w\})| \ge 2n+1-|R|- |A_2|\ge 3+N^*$.    By   minimality of $N$, $G[B \less (A_2\cup \{w\})]$ has a red $G_{i_r^*}$.  But then $G[(B \cup R)\less \{w\} ]$   has  a red $C_{2n}$, a contradiction.  \medskip

  This completes the proof of Theorem~\ref{main}. \proofsquare

  \section*{Acknowledgement}
The authors would like to thank Christian Bosse for  many helpful comments and discussion.


\begin{thebibliography}{99}

\bibitem{K5}  V.  Angeltveit,  B. D. McKay, $R(5,5)\le 48$,  to appear in J. Graph Theory. 
%
\vspace {-0.25cm}
%

\bibitem{C9C11} C. Bosse, Z-X. Song,    Multicolor Gallai-Ramsey numbers of $C_{9}$ and $C_{11}$,  submitted. arXiv:1802.06503.
%
\vspace {-0.25cm}
%
\bibitem{C13C15} C. Bosse, Z-X. Song, J. Zhang,  Improved upper bounds for Gallai-Ramsey numbers of odd cycles, submitted.
%
\vspace {-0.25cm}
%
\bibitem{DylanSong} D. Bruce, Z-X. Song, Gallai-Ramsey numbers of $C_7$ with multiple colors, submitted.
%
\vspace {-0.25cm}
%
\bibitem{CEL}  K. Cameron, J. Edmonds,   L. Lov\'asz, A note on perfect graphs, Period. Math. Hungar. 17 (1986) 173--175.
%
\vspace {-0.25cm}
%
%
%
\bibitem{chgr}
F. R. K. Chung,  R. L. Graham, Edge-colored complete graphs with precisely colored subgraphs, Combinatorica 3 (1983) 315--324.
%
%
 \vspace {-0.25cm}
%
\bibitem{R_3(10)}T. Dzido, A. Nowik, P. Szuca, New lower bound for multicolor Ramsey numbers for
even cycles, Electron.  J. Combin. 12 (2005), \#N13
%
\vspace {-0.25cm}
%

\bibitem{FLPS}
   R. J. Faudree, S. L. Lawrence, T. D. Parsons,  R. H. Schelp,  Path-cycle Ramsey numbers,
    Discrete Math. 10  (1974)   269--277.
%
\vspace {-0.25cm}
%
\bibitem{FGP}  J. Fox, A. Grinshpun,  J. Pach, The Erd\H{o}s-Hajnal conjecture for rainbow triangles,  J. Combin. Theory Ser. B 111  (2015) 75--125.
  %
\vspace {-0.25cm}
%
\bibitem{c5c6}
S. Fujita,  C. Magnant, Gallai-Ramsey numbers for cycles, Discrete Math. 311 (2011) 1247--1254.
%
\vspace {-0.25cm}
%
\bibitem{FMO}  S. Fujita, C. Magnant,  K. Ozeki, Rainbow generalizations of Ramsey theory: a survey, Graphs   Combin. 26 (2010) 1--30.
 %
\vspace {-0.25cm}
%
\bibitem{Gallai}
T. Gallai, Transitiv orientierbare Graphen, Acta Math.  Acad. Sci.  Hung. 18 (1967) 25--66.
%
%
\vspace {-0.25cm}
%
 \bibitem{GS}  A. Gy\'{a}rf\'{a}s,  G. N. S\'{a}rk\"{o}zy, Gallai colorings of non-complete graphs, Discrete Math. 310 (2010)  977--980.
%
\vspace {-0.25cm}
%
\bibitem{exponential}
A. Gy\'{a}rf\'{a}s, G.  N. S\'{a}rk\"{o}zy, A. Seb\H{o},  S. Selkow, Ramsey-type results for Gallai colorings, J. Graph Theory 64 (2010) 233--243.
%
\vspace {-0.25cm}
%
\bibitem{Gy}  A. Gy\'arf\'as, G. Simonyi,
    Edge colorings of complete graphs without tricolored triangles,
   J. Graph Theory  46   (2004) 211--216.
   %
\vspace {-0.25cm}
%
\bibitem{Hall} M. Hall, C. Magnant, K. Ozeki,  M. Tsugaki, Improved upper bounds for Gallai-Ramsey numbers of paths and cycles, J. Graph Theory 75 (2014) 59--74.
%
\vspace {-0.25cm}
%
\bibitem{KG} J. K\"orner,  G. Simonyi, Graph pairs and their entropies: modularity problems, Combinatorica 20 (2000)  227--240.
%
\vspace {-0.25cm}
%
\bibitem{K4} H. Liu, C. Magnant, A. Saito, I. Schiermeyer,  Y.  Shi,  {\small Gallai-Ramsey number for $K_4$}, submitted.
%
\vspace {-0.25cm}
%
\bibitem{GRK5} I. Schiermeyer,    personal communication.
%
\vspace {-0.25cm}
%

\bibitem{C8} Z-X. Song, J. Zhang, A conjecture on Gallai-Ramsey numbers of   even cycles and paths, submitted. arXiv:1803.07963.
%
\vspace {-0.25cm}
%
\bibitem{Rosta}  V. Rosta, On a Ramsey Type Problem of J. A. Bondy and P. Erd\H{o}s, I \& II,  J. Combin. Theory  Ser. B 15 (1973) 94--120.






\end{thebibliography}
\end{document}